\documentclass{article}
\usepackage[cp1251]{inputenc}
\usepackage[russian]{babel}
\usepackage{amssymb,amsfonts,amsmath,amsthm}
\usepackage{graphicx}
\usepackage{enumerate}
\RequirePackage{russlh}
\RequirePackage{mathlh}

\makeatletter
\newcommand*{\клей}{\nobreak\hskip\z@skip}
\renewcommand{\"}{''}
\DeclareRobustCommand*{\т}{~\textemdash{} }
\DeclareRobustCommand*{\д}{\клей\hbox{-}\клей}
\makeatother

\newdimen\theoremshape
\newdimen\defskip
\makeatletter
\renewcommand{\section}{\@startsection{section}{1}{0pt}%
{3ex plus .5ex minus .2ex}{2ex plus .2ex}%
{\center\normalfont\LARGE\bfseries\sffamily}}
\renewcommand{\subsection}{\@startsection{subsection}{2}{0pt}%
{2.75ex plus .5ex minus .2ex}{1.5ex plus .2ex}%
{\center\normalfont\large\bfseries\sffamily}}
\renewcommand{\subsubsection}{\@startsection{subsubsection}{3}{0pt}%
{2.5ex plus .5ex minus .2ex}{1ex plus .2ex}%
{\center\normalfont\bfseries\scshape}}
\makeatother
\def\postsection{.\@postskip@}
\def\postsubsection{.\@postskip@}
\def\postsubsubsection{.\@postskip@}
\def\postparagraph{.\@postskip@}
\def\postsubparagraph{.\@postskip@}

\newtheorem{theorem}{\hspace{\theoremshape} Теорема}[section]
\newtheorem{lemma}[theorem]{\hspace{\theoremshape} Лемма}
\newtheorem{prop}[theorem]{\hspace{\theoremshape} Предложение}
\newtheorem{stm}[theorem]{\hspace{\theoremshape} Утверждение}
\newtheorem{imp}{\hspace{\theoremshape} Следствие}[section]

\newenvironment{df}{\par\vskip\defskip\textbf{Определение.}}{\par\vskip\defskip}

\newenvironment{denote}{\par\vskip\defskip\textbf{Обозначение.}}{\par\vskip\defskip}

\renewenvironment{proof}{\par $\square\quad$}{$\blacksquare$ \par\vskip\defskip}

\newenvironment{nums}[1]
{\begin{enumerate}
\setlength\itemsep{#1pt}}
{\end{enumerate}}
\newcommand{\eqn}[1]{\begin{equation}#1\end{equation}}
\newcommand{\equ}[1]{\begin{equation*}#1\end{equation*}}
\newcommand{\case}[1]{\begin{cases}#1\end{cases}}
\makeatletter
\@addtoreset{equation}{section}
\makeatother

\renewcommand{\ge}{\geqslant}
\renewcommand{\le}{\leqslant}
\newcommand{\fa}{\,\forall\,}
\newcommand{\bes}{\infty}
\newcommand{\es}{\varnothing}
\newcommand{\subs}{\subset}
\newcommand{\sups}{\supset}
\newcommand{\subsneq}{\subsetneqq}
\newcommand{\sm}{\setminus}
\newcommand{\cln}{\colon}
\newcommand{\nl}{\lhd}
\newcommand{\wg}{\wedge}
\newcommand{\lhdp}{\leftthreetimes}
\newcommand{\Lra}{\Leftrightarrow}
\newcommand{\xra}{\xrightarrow}
\newcommand{\us}{\underset}
\newcommand{\ol}{\overline}
\newcommand{\wt}{\widetilde}

\newcommand{\suml}[2]{\sum\limits_{{#1}}^{{#2}}}
\newcommand{\sums}[1]{\sum\limits_{{#1}}}
\newcommand{\sumiun}{\sum\limits_{i=1}^{n}}
\newcommand{\prods}[1]{\prod\limits_{{#1}}}
\newcommand{\oplusl}[2]{\bigoplus\limits_{{#1}}^{{#2}}}
\newcommand{\opluss}[1]{\bigoplus\limits_{{#1}}^{}}
\newcommand{\caps}[1]{\bigcap\limits_{{#1}}}

\newcommand*{\bw}[1]{#1\nobreak\discretionary{}{\hbox{$\mathsurround=0pt #1$}}{}}
\newcommand{\sco}{,\ldots,}
\newcommand{\spl}{\bw+\ldots\bw+}
\newcommand{\sop}{\bw\oplus\ldots\bw\oplus}
\newcommand{\sd}{\bw\cdot\ldots\bw\cdot}

\newcommand{\ha}[1]{\left\langle#1\right\rangle}
\newcommand{\ba}[1]{\bigl\langle#1\bigr\rangle}
\newcommand{\hr}[1]{\left(#1\right)}
\newcommand{\br}[1]{\bigl(#1\bigr)}
\newcommand{\Br}[1]{\Bigl(#1\Bigr)}

\newcommand{\bm}[1]{\bigl|#1\bigr|}
\newcommand{\hn}[1]{\left\|#1\right\|}
\newcommand{\bn}[1]{\bigl\|#1\bigr\|}
\newcommand{\hs}[1]{\left[#1\right]}

\newcommand{\hc}[1]{\left\{#1\right\}}
\newcommand{\bc}[1]{\bigl\{#1\bigr\}}

\newcommand{\mbb}{\mathbb}
\newcommand{\mbf}{\mathbf}
\newcommand{\mcl}{\mathcal}
\newcommand{\mfr}{\mathfrak}
\newcommand{\R}{\mbb{R}}

\newcommand{\Z}{\mbb{Z}}
\newcommand{\N}{\mbb{N}}
\newcommand{\T}{\mbb{T}}
\newcommand{\F}{\mbb{F}}
\newcommand{\Cbb}{\mbb{C}}
\newcommand{\Hbb}{\mbb{H}}
\newcommand{\idb}{\mbf{1}}
\newcommand{\ib}{\mbf{i}}
\newcommand{\jb}{\mbf{j}}
\newcommand{\kb}{\mbf{k}}
\newcommand{\Zc}{\mcl{Z}}
\newcommand{\ggt}{\mfr{g}}
\newcommand{\hgt}{\mfr{h}}

\newcommand{\al}{\alpha}
\newcommand{\ga}{\gamma}
\newcommand{\Ga}{\Gamma}
\newcommand{\de}{\delta}
\newcommand{\De}{\Delta}
\newcommand{\ep}{\varepsilon}
\newcommand{\la}{\lambda}
\newcommand{\La}{\Lambda}
\newcommand{\rh}{\rho}
\newcommand{\si}{\sigma}
\newcommand{\ta}{\theta}
\newcommand{\ph}{\varphi}
\newcommand{\om}{\omega}
\newcommand{\Om}{\Omega}

\DeclareMathOperator{\Lie}{Lie}
\DeclareMathOperator{\Ker}{Ker}
\DeclareMathOperator{\Ad}{Ad}
\DeclareMathOperator{\rk}{rk}
\DeclareMathOperator{\diag}{diag}
\DeclareMathOperator{\id}{id}

\newcommand{\GL}{\mbf{GL}}
\newcommand{\SL}{\mbf{SL}}
\newcommand{\Or}{\mbf{O}}

\newcommand{\SO}{\mbf{SO}}

\begin{document}

\author{О. Г. Стырт}
\title{<<О пространстве орбит компактной линейной группы Ли\\
с коммутативной связной компонентой>>}
\date{}

\maketitle

\tableofcontents

\section{Введение}\label{introd}

\begin{df} Непрерывное отображение гладких многообразий назовём \textit{кусочно-гладким}, если оно переводит любое гладкое подмногообразие в конечное объединение гладких подмногообразий.
\end{df}

В частности, всякое собственное гладкое отображение гладких многообразий является кусочно-гладким.

Рассмотрим дифференцируемое действие некоторой компактной группы Ли $G$ на гладком многообразии~$M$.

\begin{df} Будем говорить, что фактор действия $G\cln M$ \textit{диффеоморфен} (\textit{кусочно-диффеоморфен}) гладкому многообразию~$M'$, если топологический фактор $M/G$ гомеоморфен~$M'$, причём отображение факторизации $M\to M'$ гладкое (кусочно-гладкое).
\end{df}

\begin{df} Будем говорить, что фактор действия $G\cln M$ является \textit{гладким многообразием}, если он кусочно-диффеоморфен некоторому гладкому
многообразию.
\end{df}

Пусть имеется точное линейное представление компактной группы Ли $G$ в вещественном векторном пространстве~$V$. Нас будет интересовать вопрос о том, является ли фактор $V/G$ этого действия топологическим многообразием, а также является ли он гладким многообразием. Для краткости будем в дальнейшем называть топологическое многообразие просто <<многообразием>>.

\begin{denote} Для произвольного линейного представления $G\cln V$ и целого неотрицательного числа $d$ под записью $(V\oplus\R^d)/G$ будем, если не оговорено противное, подразумевать фактор линейного представления
\eqn{\label{act}
G\cln V\oplus\R^d,\;g\cln v+x\to gv+x,\;g\in G,\,v\in V,\,x\in\R^d.}
\end{denote}

\begin{df} Линейный оператор в пространстве над некоторым полем называется \textit{отражением} (соотв. \textit{псевдоотражением}), если подпространство его неподвижных точек имеет коразмерность $1$~(соотв. $2$).
\end{df}

Для случая, когда группа $G$ конечна, М. А. Михайловой был получен следующий результат~(\cite{MAMich}): если фактор $V/G$ диффеоморфен пространству $V$, то группа $G$ порождена псевдоотражениями, а если $G$ порождается псевдоотражениями, то фактор $V/G$ гомеоморфен $V$, в частности, является многообразием. Сформулируем обобщение этого результата, которое будет доказано в п.~\ref{facts}.

\begin{theorem}\label{Mich} Пусть группа $G$ конечна. Тогда если $(V\oplus\R^d)/G$\т гладкое многообразие для некоторого $d\ge0$, то группа $G$ порождается псевдоотражениями, а если $G$ порождена псевдоотражениями, то $V/G$\т многообразие.
\end{theorem}

Через $G^0$ будем обозначать связную компоненту единицы группы~$G$, а через $\ggt$\т её касательную алгебру. В пространстве $V$ можно зафиксировать $G$\д инвариантное скалярное умножение. Тогда группа $G$ действует ортогональными операторами: $G\subs\Or(V)$, причём $G^0\subs\SO(V)$.

Пусть $G_v$\т стабилизатор вектора $v\in V$, а $\ggt_v:=\Lie G_v=\{\xi\in\ggt\cln\xi v=0\}$. Группа $G_v$ переводит в себя вектор $v$, орбиту $Gv$, а значит, и касательное пространство $T_v(Gv)=\ggt v$, а также подпространство $N_v:=(\ggt v)^{\perp}$. Через $M_v$ будем обозначать ортогональное дополнение в $N_v$ к подпространству $N_v^{G_v}$ неподвижных векторов для действия~$G_v\cln N_v$. Тогда $V=\ggt v\oplus N_v^{G_v}\oplus M_v$, $G_vM_v=M_v$, а
\eqn{\label{NvGv}
N_v/G_v\cong N_v^{G_v}\times(M_v/G_v).}
Пространство $M_v$ не содержит ненулевых векторов, неподвижных относительно $G_v$; в частности, если оно нетривиально, то группа $G_v$ не может действовать на нём тождественно.

Для произвольного элемента $g\in G$ введём обозначение
\equ{\label{ome}
\om(g):=\rk(E-g)-\rk\br{E-\Ad(g)}.}
Положим
\equ{\label{omega}
\Om:=\bc{g\in G\cln\om(g)\in\{0;2\}}\subs G.}

Пусть $v\in V$\т некоторый вектор с конечным стабилизатором. Тогда $\ggt_v=0$, и для любого $g\in G_v$ имеем
\eqn{\label{dims}
\dim\br{(E-g)N_v}=\dim\br{(E-g)V}-\dim\br{(E-g)(\ggt v)}.}
Отображение $\ggt\to\ggt v,\,\xi\to\xi v$ является линейным изоморфизмом пространств $\ggt$ и $\ggt v$, и если с его помощью эти пространства отождествить, то элемент $g$ действует на $\ggt v$ так же, как оператор $\Ad(g)$\т на~$\ggt$. В таком случае
$\dim\br{(E-g)(\ggt v)}=\rk\br{E-\Ad(g)}$, и, в силу~\eqref{dims}, $\dim\br{(E-g)N_v}=\om(g)$. В частности, элемент $g\in G_v$ принадлежит $\Om$, если и только если он действует на $N_v$ псевдоотражением или тождественно.

Для подпространства $W\subs V$ группа $\Or(W)$ вкладывается в группу $\Or(V)$ как подгруппа всех ортогональных преобразований пространства $V$, действующих тождественно на~$W^{\perp}$. Если $g\in\Or(V)$, то $\Or(gW)=g\Or(W)g^{-1}$. Далее, пусть $G[W]:=G\cap\Or(W)$\т подгруппа Ли всех элементов из $G$, действующих тождественно на~$W^{\perp}$. Алгебра $\ggt[W]:=\Lie G[W]\subs\ggt$ состоит из всех элементов алгебры $\ggt$, переводящих
$W^{\perp}$ в нуль. Для всякого $g\in G$ имеем $G[gW]=gG[W]g^{-1}$ и $\ggt[gW]=\Ad(g)\ggt[W]$.

В данной работе будет рассмотрена ситуация, когда группа $G^0$ коммутативна, то есть является тором. Очевидно, что это свойство сохраняется при переходе
к факторгруппе.

Поскольку $G^0\nl G$, всякое $G^0$\д инвариантное подпространство под действием любого $g\in G$ переходит вновь в $G^0$\д инвариантное; неприводимый $G^0$\д подмодуль\т в неприводимый. Кроме того, образы изоморфных $G^0$\д подмодулей изоморфны. Пространство $V$ можно разложить в прямую сумму попарно ортогональных неприводимых $G^0$\д инвариантных подпространств. При этом прямые суммы изоморфных неприводимых подмодулей (изотипные компоненты) также попарно ортогональны, а группа $G$ переставляет их.

Все неприводимые представления коммутативной группы $G^0$ имеют размерность $1$ или~$2$. Всякое двумерное представление имеет
комплексную структуру, в соответствии с которой группа $G^0$ действует на комплексной прямой умножением на числа из группы Ли $\T:=\{z\in\Cbb\cln|z|=1\}$. Алгебру $\Lie\T$ можно отождествить с~$\R$. Каждому двумерному неприводимому представлению группы $G^0$ соответствует вес $\la\cln G^0\to\T$\т нетривиальный гомоморфизм групп Ли, а также его дифференциал\т нетривиальная линейная функция $\la\cln\ggt\to\R$. Пространство $\ggt$ можно отождествить с $\ggt^*$ при помощи скалярного умножения, инвариантного относительно $\Ad(G)$, поэтому линейную функцию $\la$ можно понимать как вектор пространства~$\ggt$. Вес $\la\cln G^0\to\T$ (вес $\la\in\ggt$) определён с точностью до поточечного обращения (смены знака), так
как на комплексной прямой можно сделать замену координаты $z\to\ol{z}$. Одномерным представлениям группы $G^0$ естественно сопоставить тривиальные веса $\la\cln G^0\to\T$ и нулевые веса $\la=0\in\ggt$, поскольку на этих представлениях группа $G^0$ действует тождественно.

Классы изоморфных неприводимых представлений группы $G^0$ характеризуются весами $\la\in\ggt$, определёнными с точностью до знака. Множество $P$ весов $\la$, соответствующее разложению $V$ в прямую сумму неприводимых компонент (с учётом кратностей последних), не зависит от разложения и переходит в себя с точностью до знака (но с учётом кратностей) под действием любого оператора $\Ad(g)$,~$g\in G$.

Пусть $V_{\la}$ ($\la\in P$)\т изотипная компонента, соответствующая неприводимым представлениям с весами~$\pm\la$. В частности, $V_0$\т изотипная компонента, соответствующая одномерным неприводимым представлениям. Она переходит в себя под действием $G$, причём под действием $G^0$\т поточечно. Имеем $V_{\la_1}=V_{\la_2}\Lra\la_1=\pm\la_2$. Для подмножества $Q\subs P$ обозначим через $V_Q$ сумму изотипных компонент $V_{\la}$
по всем~$\la\in Q$.

Всякая изотипная компонента $V_{\la}$, $\la\ne0$, имеет комплексную структуру, причём элемент $g\in G^0$ действует на ней умножением на $\la(g)$, а $\xi\in\ggt$\т на~$i\la(\xi)$. Под действием элемента $g\in G$ компонента $V_{\la}$ переходит в компоненту $V_{\Ad(g)\la}$, а подпространство $V_Q$\т в~$V_{\Ad(g)Q}$. При этом компонента $V_{\la}$ ($\la\ne0$) переходит в себя тогда и только тогда, когда $\Ad(g)\la=\pm\la$. Если в последнем равенстве стоит знак <<плюс>> (<<минус>>), то ограничение $g$ на $V_{\la}$ коммутирует (антикоммутирует) с каждым $\xi\in\ggt$ и потому линейно (антилинейно) над полем~$\Cbb$. Поскольку представление точное, пересечение ядер всех весов $G^0\to\T$ из $P$ тривиально, а пересечение ядер соответствующих весов $\la\in\ggt$ нулевое, и, значит, эти веса линейно порождают $\ggt$, то есть~$\ha{P}=\ggt$.

Для $G^0$\д инвариантного подпространства $W\subs V$ группа $G[W]\cap G^0$ (алгебра $\ggt[W]$) есть пересечение ядер весов $G^0\to\T$ ($\ggt\to\R$), соответствующих~$W^{\perp}$.

Разложением множества векторов конечномерного линейного пространства на компоненты будем называть его представление в виде объединения своих подмножеств, линейные оболочки которых линейно независимы. Если среди этих линейных оболочек по крайней мере две нетривиальны, то такое разложение назовём собственным. Будем говорить, что множество \textit{неразложимо}, если оно не допускает ни одного собственного разложения на компоненты. Всякое множество векторов разлагается на неразложимые компоненты единственным образом (с точностью до распределения нулевого вектора), причём для любого его
разложения на компоненты каждая компонента является объединением некоторых его неразложимых компонент (вновь с точностью до нулевого вектора).

\begin{df}\label{stay} Конечное множество векторов конечномерного пространства, рассматриваемое с учётом кратностей своих элементов, назовём \textit{$q$\д устойчивым} ($q\ge0$), если его линейная оболочка не меняется при удалении из него любых векторов в количестве не более $q$ (с учётом кратностей).
\end{df}

Добавление и удаление нулевых векторов не влияет на $q$\д устойчивость множества. Не влияет на неё и умножение некоторых векторов множества на ненулевые скаляры. Ясно, что образ $q$\д устойчивого множества при линейном отображении $q$\д устойчив. В частности, если конечное множество линейных функций на некотором пространстве $q$\д устойчиво, то множество их ограничений на любое подпространство также $q$\д устойчиво.

Легко видеть, что множество, содержащее не более $q$ ненулевых векторов, $q$\д устойчиво тогда и только тогда, когда оно не содержит ненулевых векторов. Значит, любое $q$\д устойчивое множество с $k$\д мерной линейной оболочкой ($k>0$) содержит более $q$ ненулевых векторов, и при удалении $q$ из них линейная оболочка должна остаться прежней, $k$\д мерной. Следовательно, в таком множестве не может быть менее $k+q$ ненулевых векторов.

Если $q$\д устойчивое множество разложить на компоненты (не обязательно неразложимые), то каждая из компонент $q$\д устойчива как образ исходного множества при линейном отображении\т проекции на соответствующее прямое слагаемое (с точностью до добавления и удаления нулей).

Перейдём к формулировке результатов, которые будут получены в дальнейшем. Прежде всего, множество всех ненулевых весов в $P$ можно разложить в дизъюнктное объединение неразложимых компонент~$\La\subs P$. Пусть $L$\т множество подпространств, состоящее из всех $V_{\La}$ ($\La$\т неразложимая компонента), а также~$V_0$. Тогда все подпространства из $L$ попарно ортогональны, а $V$ есть их прямая сумма.

Пространство $\ggt$ является прямой суммой линейных оболочек всех неразложимых компонент~$\La\subs P$. Произведение всех подгрупп $G[W]$, $W\in L$, прямое и содержит $G^0$, так как прямая сумма всех подпространств $\ggt\hs{V_{\La}}$ есть~$\ggt$. Точному представлению $G\hs{V_{\La}}\cln V_{\La}$ соответствует неразложимое множество весов.

В п.~\ref{2stab} будет доказана следующая теорема.

\begin{theorem}\label{submain} Допустим, что $P$\т $2$\д устойчивое множество, а $V/G$\т гладкое многообразие. Тогда группа $G$ совпадает с прямым произведением всех подгрупп $G[W]$, $W\in L$. Кроме того, для любого $W\in L$ найдётся число $d\ge0$, такое что $(W\oplus\R^d)/G[W]$\т гладкое многообразие.
\end{theorem}

Если группа $G$ совпадает с прямым произведением всех подгрупп $G[W]$, $W\in L$, а $W/G[W]$\т многообразие для каждого $W\in L$, то и фактор $V/G$ является многообразием как декартово произведение факторов $W/G[W]$,~$W\in L$. (Это верно и без требования $2$\д устойчивости~$P$.)

Если $P$\т $2$\д устойчивое множество, то каждая компонента $\La$ сама $2$\д устойчива, а множество весов, соответствующих действию
$G\hs{V_{\La}}\cln V_{\La}$, неразложимо, $2$\д устойчиво и не содержит нулей.

В п.~\ref{simpl} вопрос о том, является ли (гладким) многообразием фактор произвольного представления группы Ли рассматриваемого типа, будет сведён к аналогичному вопросу для представления с $2$\д устойчивым множеством весов.

В случае $2$\д устойчивого множества $P$ нужно разложить $V$ в прямую сумму подпространств $W\in L$, проверить, совпадает ли $G$ с прямым произведением подгрупп $G[W]$, $W\in L$ (если это не так, то $V/G$ не есть гладкое многообразие), и изучить отдельно каждое действие
$G[W]\cln W$,~$W\in L$. При $W=V_0$ можно воспользоваться теоремой~\ref{Mich}, так как группа $G[W]$ конечна.

Теперь требуется разобраться с ситуацией $2$\д устойчивого неразложимого множества $P$, не содержащего нулей. Следующая теорема описывает этот случай при дополнительном условии~$\dim G>1$.

\begin{theorem}\label{main} Предположим, что множество весов $P$ неразложимо, $2$\д устойчиво и не содержит нулей.
\begin{nums}{-1}
\item Если $\dim G>1$, а $(V\oplus\R^d)/G$\т гладкое многообразие для некоторого $d\ge0$, то выполнены четыре условия:
\begin{nums}{-1}
\item\label{m2stab} множество $P$ является $2$\д устойчивым и содержит $m+2$ ненулевых веса ($m=\dim G>0$);
\item\label{oplu} пространство $V$ разлагается в прямую сумму попарно ортогональных двумерных $G^0$\д инвариантных подпространств $W_1\sco W_{m+2}$,
переставляемых группой $G$, причём последняя переводит в себя подпространство $W_1\sop W_m$, а при $m>2$\т все $W_i$, $i=1\sco m+2$;
\item\label{adj} найдётся элемент $g\in G$, переводящий в себя все $W_i$, для которого $\Ad(g)=-E$;
\item\label{finst} для любого вектора $v\in V$ с конечным стабилизатором $G_v=\ha{G_v\cap\Om}$.
\end{nums}
\item Если $\dim G>0$ (не обязательно $\dim G>1$) и условия~\ref{m2stab}~---~\ref{finst} выполнены, то $V/G$\т многообразие.
\end{nums}
\end{theorem}

\begin{theorem}\label{mainp} Если для представления $G\cln V$ выполнены условия~\ref{m2stab}~---~\ref{adj} теоремы~\ref{main} и при этом группа $G$ переводит в себя все подпространства $W_1\sco W_{m+2}$, то $V/G$\т многообразие.
\end{theorem}

Пусть теперь множество $P$ неразложимо, $2$\д устойчиво и не содержит нулей, а группа $G$ одномерна (или, что равносильно, $\dim G=1$, а $P$ содержит по крайней мере $3$ ненулевых веса и ни одного нулевого). Пространство $V$ имеет комплексную структуру, а $\Ad(G)\subs\{\pm E\}$. Элемент $g\in G$, для которого $\Ad(g)=E$ ($\Ad(g)=-E$), действует на $V$ линейно (антилинейно) над~$\Cbb$. В п.~\ref{examp} будет доказано, что все комплексные отражения группы $G$ порождают в ней конечную подгруппу $H\nl G$, фактор $V/H$ диффеоморфен $V$, и на нём группа $G/H$ действует линейно. Очевидно,
что все операторы группы $H$ комплексно-линейны.

В п.~\ref{1dimg} мы разберём указанный случай и докажем для него теоремы~\ref{main1}~---~\ref{GiHr}.

\begin{theorem}\label{main1} Если $(V\oplus\R^d)/G$\т гладкое многообразие для некоторого $d$, то $\dim_{\Cbb}V=3$, $\Ad(G)=\{\pm E\}$, группа $G$ порождена множеством $\Om$, а представление $(G/H)\cln(V/H)$ приводимо. При выполнении всех перечисленных условий фактор $V/G$ является многообразием.
\end{theorem}

При формулировке остальных трёх теорем будем для удобства считать, что $\dim_{\Cbb}V=3$, $\Ad(G)=\{\pm E\}$.

\begin{theorem}\label{GrHi} Предположим, что либо группа $H$ неприводима над полем $\Cbb$, либо $G$ приводима над~$\R$. Тогда если $(V\oplus\R^d)/G$\т гладкое многообразие для некоторого $d$, то $G=\ha{\Om}$, а если $G=\ha{\Om}$, то $V/G$\т многообразие.
\end{theorem}

Следуя~\cite{ShT}, для произвольных $p,q\in\N$ обозначим через $G(pq,p,3)$ конечную подгруппу в $\GL_3(\Cbb)$, порождённую всеми операторами перестановок координат, а также всеми операторами видов $\diag(\ep_1,1,1)$ и $\diag(1,\ep_2,\ep_2^{-1})$, $\ep_1^q=\ep_2^{pq}=1$. Группа $G(pq,p,3)$ переставляет координатные прямые, а все её операторы, осуществляющие на множестве этих прямых чётную подстановку, образуют в ней подгруппу $G'(pq,p,3)$ индекса~$2$. Пусть $G'(q,3)\subs\GL_{\Cbb}(V)$, $q\in\N$\т одномерная группа Ли, порождённая группой $G'(q,1,3)$ и всеми скалярными операторами $\la E$,~$\la\in\T$. В её нормализаторе лежит инволютивный оператор $g_{(2,3)}\cln(z_1,z_2,z_3)\to(\ol{z_1},\ol{z_3},\ol{z_2})$,
антилинейный над~$\Cbb$. Через $G(q,3)\subs\GL_{\R}(V)$ обозначим одномерную группу Ли, порождённую $G'(q,3)$ и~$g_{(2,3)}$.

\begin{theorem}\label{HG3} Если группа $H$ неприводима над полем $\Cbb$ и не совпадает ни с одной из групп $G(3q,3,3)$, $q\in\N$, то $V/G$\т многообразие.
\end{theorem}

\begin{theorem}\label{GiHr} Предположим, что группа $G$ неприводима над $\R$, а $H$\т приводима над~$\Cbb$. Тогда если для некоторого $d$ фактор
$(V\oplus\R^d)/G$ является гладким многообразием, то $G=G(q,3)$, $q>1$, а если $G=G(q,3)$, $q>1$, то $V/G$\т многообразие.
\end{theorem}

В тех случаях, когда утверждается, что $V/G$\т многообразие, вопрос о том, является ли этот фактор гладким многообразием, остаётся открытым.

\section{Вспомогательные факты}\label{facts}

\begin{stm} Пусть $M_1$, $M_2$, $M_3$\т гладкие многообразия, $f_1\cln M_1\to M_2$ и $f_2\cln M_2\to M_3$\т непрерывные отображения, причём отображение $f_1$ является собственным, гладким и сюръективным. Тогда отображения $f_2$ и $f:=f_2\circ f_1\cln M_1\to M_3$ являются или не являются кусочно-гладкими одновременно.
\end{stm}

\begin{proof} Образ (прообраз) под действием $f_1$ любого гладкого подмногообразия из $M_1$ (из $M_2$) является конечным объединением гладких
подмногообразий в $M_2$ (в~$M_1$). Кроме того, $f(U_1)=f_2\br{f_1(U_1)}$ и $f_2(U_2)=f\br{f_1^{-1}(U_2)}$ для любых подмножеств~$U_i\subs M_i$.
\end{proof}

\begin{imp}\label{redu} Пусть имеется линейное представление компактной группы Ли $G$ (не обязательно с коммутативной связной компонентой) в векторном пространстве~$V$. Предположим, что для некоторой нормальной подгруппы $H\nl G$ фактор $V/H$ диффеоморфен гладкому многообразию~$M$. Тогда факторы $M/G$ и $V/G$ являются или не являются гладкими многообразиями одновременно.
\end{imp}

Вновь рассмотрим линейное представление $G\cln V$ (по-прежнему не требуем коммутативность~$G^0$). Отображение факторизации $V\to V/G$ обозначим через~$\pi$.

\begin{prop} Если $V/G$\т гладкое многообразие, то в $V$ найдётся связная $G$\д инвариантная окрестность нуля $U$, для которой фактор $U/G$ кусочно-диффеоморфен открытому шару размерности~$\dim(V/G)$.
\end{prop}

\begin{proof} У точки $\pi(0)\in V/G$ есть окрестность, диффеоморфная клетке, то есть в $V$ существует $G$\д инвариантная окрестность нуля $U$, такая что фактор $U/G$ кусочно-диффеоморфен открытому шару размерности $\dim(V/G)$ и, в частности, связен. В окрестности $U$ всякая компонента линейной связности открыта. Если $U^0$\т компонента линейной связности, содержащая $0$, то подмножества $U^0$ и $U\sm U^0$ открыты и $G$\д инвариантны. Их
образы при факторизации открыты, не пересекаются и дают в объединении связный фактор $U/G$. Следовательно, окрестность $U$ совпадает с $U^0$, а значит, связна.
\end{proof}

\begin{lemma}\label{locman} Если некоторая окрестность точки $\pi(0)$ в $V/G$ является (гладким) многообразием, то $V/G$\т (гладкое) многообразие.
\end{lemma}

\begin{proof} В пространстве $V$ найдётся $G$\д инвариантная окрестность нуля $U$, для которой фактор $U/G$\т (гладкое) многообразие; множество $U$ содержит некоторый открытый шар $B$ с центром в нуле. Этот шар $G$\д инвариантен, и $B/G$\т (гладкое) многообразие. Поскольку существует диффеоморфизм между $V$ и $B$, перестановочный с действием $G$, фактор $V/G$ также является (гладким) многообразием.
\end{proof}

\begin{theorem}\label{slice} Пусть $v\in V$\т некоторый вектор. Фактор $V/G$ является многообразием локально в точке $\pi(v)$ тогда и только тогда, когда $N_v/G_v$\т многообразие. Кроме того, если $V/G$\т гладкое многообразие, то и $N_v/G_v$\т гладкое многообразие.
\end{theorem}

\begin{proof} Положим~$H:=G_v$. По теореме о слайсе~(\cite[~гл.~II,~~\S~4\т 5]{Bredon}) в пространстве $N_v$ существует $H$\д инвариантная окрестность $U$ точки $v$, такая что подмножество $GU$ открыто в $V$ и является однородным расслоением $G\us{H}*U$, то есть
\eqn{\label{hoge}
U\cap gU=\es\;\fa g\in G\sm H.}
Если $V/G$\т гладкое многообразие, то и $(GU)/G$\т гладкое многообразие, так как $GU$\т открытое подмножество в~$V$. Отображение факторизации $GU\to(GU)/G$ при ограничении на $U$ даёт факторизацию $U$ по~$H$. Значит, и фактор $U/H$ является гладким многообразием. Но тогда в пространстве $N_v$ имеется $H$\д инвариантная окрестность нуля $U-v$, фактор которой по действию $H$ является гладким многообразием: отображение $U\to U-v$,
$w\to w-v$ является диффеоморфизмом, перестановочным с действием~$H$.

Из теоремы о слайсе также следует, что фактор $V/G$ локально в точке $\pi(v)$ гомеоморфен фактору $N_v/G_v$ локально в нуле.

Осталось применить лемму~\ref{locman} и получить оба утверждения теоремы.
\end{proof}

\begin{imp}\label{slicd} Пусть $v\in V$\т произвольный вектор. Тогда если $(V\oplus\R^d)/G$\т гладкое многообразие при некотором $d\ge0$, то и $(N_v\oplus\R^d)/G_v$\т гладкое многообразие.
\end{imp}

\begin{proof} При естественном вложении $V$ в пространство $V\oplus\R^d$, на котором задано действие~\eqref{act}, ортогональным дополнением к
подпространству $\ggt v$ будет в точности~$N_v\oplus\R^d$. Осталось воспользоваться теоремой~\ref{slice}.
\end{proof}

\begin{imp}\label{trans} Если для некоторого вектора $v\in V$ подпространство $M_v$ нетривиально, а группа $G_v$ действует транзитивно на единичной сфере в нём, то $V/G$\т не многообразие.
\end{imp}

\begin{proof} Ясно, что фактор пространства $M_v$ по действию $G_v$ гомеоморфен полупрямой $\R_{\ge0}$, а фактор $N_v/G_v$, согласно~\eqref{NvGv}, гомеоморфен $N_v^{G_v}\times\R_{\ge0}$\т замкнутому полупространству\т и, таким образом, не является многообразием. Из теоремы~\ref{slice} следует,
что и $V/G$\т не многообразие.
\end{proof}

В пространстве $V$ можно выделить \textit{страты}\т подмножества векторов со стабилизатором, сопряжённым фиксированной подгруппе в~$G$. Всего (непустых) стратов конечное число, каждый из них $G$\д инвариантен, $V$\т их дизъюнктное объединение. В соответствии с этим фактор $V/G$ также распадается в дизъюнктное объединение образов стратов при факторизации, которые мы тоже будем называть стратами в~$V/G$. Для любого вектора $v\in V$ страт в $V$, содержащий $v$,\т гладкое подмногообразие в $V$ коразмерности~$\dim M_v$. Образ этого страта при факторизации также является гладким многообразием и имеет коразмерность в $V/G$, равную размерности фактора~$M_v/G_v$. В частности, коразмерность страта в $V/G$ не превосходит коразмерности соответствующего страта в $V$, поскольку $\dim(M_v/G_v)\le\dim M_v$. Ограничение отображения $\pi$ на любой страт является локально
тривиальным расслоением. Страт коразмерности $0$ (то есть открытый) в пространстве $V$ имеется ровно один, ему соответствует единственный открытый страт в~$V/G$. Эти страты будем называть \textit{главными}.

\begin{lemma}\label{strat1} Если фактор $V/G$ является многообразием, то в нём нет стратов коразмерности~$1$.
\end{lemma}

\begin{proof} Предположим, что найдётся страт коразмерности $1$, и пусть $v$\т содержащийся в нём вектор. Тогда $\dim(M_v/G_v)=1$, значит, пространство $M_v$ нетривиально, а фактор его единичной сферы по действию $G_v$ дискретен. В силу следствия~\ref{trans}, группа $G_v$ не может транзитивно действовать на единичной сфере в~$M_v$. Таким образом, фактор указанного действия состоит более чем из одной точки и, будучи дискретным, не является связным. Последнее имеет место лишь в том случае, когда пространство $M_v$ одномерно, а группа $G_v$ действует на нём тождественно, что, в свою очередь, невозможно.
\end{proof}

\begin{stm}\label{MvW} В любом $G^0$\д инвариантном подпространстве $W\subs V$ существует вектор $v$, для которого $M_v\subs W^{\perp}$.
\end{stm}

\begin{proof} Среди векторов пространства $W$ выберем вектор $v$ с минимальным по включению стабилизатором. В пространстве $N_v$ существует окрестность $U$ вектора $v$, удовлетворяющая~\eqref{hoge}. В частности, стабилизатор любого вектора из $U$ содержится в~$G_v$. Следовательно,
стабилизатор любого вектора из $W\cap U$ совпадает с $G_v$, $v\in W\cap U\subs N_v^{G_v}$, $W\cap N_v=\ha{W\cap U}\subs N_v^{G_v}$,
$W=\ggt v\oplus(W\cap N_v)\subs\ggt v\oplus N_v^{G_v}=M_v^{\perp}$.
\end{proof}

\begin{lemma}\label{transfer} Пусть $W\subs V$\т ненулевое $G^0$\д инвариантное подпространство, причём на его единичной сфере группа $G[W]$ действует транзитивно. Тогда $V/G$\т не многообразие.
\end{lemma}

\begin{proof} Выберем вектор $v\in W^{\perp}$, для которого $M_v\subs W$. Согласно условию, в $W$ нет ненулевых векторов, неподвижных относительно~$G[W]$. Кроме того, $G[W]\subs G_v$, откуда $W\cap N_v^{G_v}=0$. При этом $M_v\subs W\subs N_v$, значит, $W=M_v$. Из условия также следует, что действие группы $G_v$, содержащей $G[W]$, на единичной сфере в ненулевом подпространстве $M_v$ транзитивно. Остаётся применить следствие~\ref{trans} и получить требуемое.
\end{proof}

\begin{imp}\label{1dim} Если $V/G$\т многообразие, то в $G$ нет отражений относительно $G^0$\д инвариантных гиперплоскостей.
\end{imp}

\begin{proof} Если $g\in G$\т отражение относительно $G^0$\д инвариантной гиперплоскости $W'$, то на $G^0$\д инвариантной прямой $W:=(W')^{\perp}$ оператор $g\in G[W]$ действует умножением на~$(-1)$. Поэтому группа $G[W]$ транзитивно действует на единичной сфере в пространстве $W$ (состоящей из двух точек).
\end{proof}

\begin{imp}\label{2dim} Если $V/G$\т многообразие, то для всякого двумерного $G^0$\д инвариантного подпространства $W$ группа $G[W]$ конечна.
\end{imp}

\begin{proof} Группа $G[W]$ вкладывается в $\Or(W)\cong\Or_2$. Если она при этом бесконечна, то содержит группу $\SO_2$ всех поворотов двумерного пространства $W$ и транзитивно действует на его единичной сфере (окружности).
\end{proof}

Через $U_k$ ($k>0$) обозначим объединение всех стратов в $V$, образы которых под действием $\pi$ имеют коразмерность более $k$ в факторе~$V/G$. В таком случае $U_k/G$\т объединение всех стратов в $V/G$ коразмерности более~$k$, а $U_k$\т объединение некоторых стратов в $V$, коразмерность каждого из которых также более~$k$.

Если $V/G$\т гладкое многообразие, то существует связная $G$\д инвариантная окрестность нуля $U\subs V$, такая что фактор $U/G$ кусочно-диффеоморфен открытому шару размерности~$\dim(V/G)$. В частности, этот фактор является связным гладким многообразием, все его гомотопические группы тривиальны, а его подмножество $(U\cap U_k)/G$ есть объединение конечного числа гладких подмногообразий коразмерностей более~$k$. Следовательно, $(U\sm U_k)/G$\т связное топологическое пространство с тривиальными гомотопическими группами~$\pi_1\sco\pi_{k-1}$.

В связном гладком многообразии $U\subs V$ подмножество $U\cap U_k$ представляет собой объединение конечного числа гладких подмногообразий коразмерностей более $k$, поэтому топологическое пространство $U\sm U_k$ непусто и связно, отсюда фактор $(U\sm U_k)/G^0$ также непуст и связен. Кроме того, для любого открытого шара $B\subs V$ подмножество $B\sm U_k$ в пространстве $V$ непусто, связно и имеет тривиальные гомотопические
группы~$\pi_1\sco\pi_{k-1}$.

\begin{theorem}\label{gener} Если $V/G$\т гладкое многообразие, то факторгруппа $G/G^0$ порождается смежными классами, содержащими представителя, оставляющего на месте некоторый вектор, которому соответствует страт в $V/G$ коразмерности не более $2$,\т то есть вектор из~$V\sm U_2$.
\end{theorem}

\begin{proof} Рассмотрим связную $G$\д инвариантную окрестность нуля  $U\subs V$, такую что фактор $U/G$ кусочно-диффеоморфен открытому шару размерности~$\dim(V/G)$. Как уже было выяснено, топологическое пространство $(U\sm U_2)/G^0$ непусто и связно, а его фактор $(U\sm U_2)/G$ по действию конечной группы $G/G^0$ односвязен. Поэтому факторгруппа $G/G^0$ порождена элементами, имеющими неподвижную точку в~$(U\sm U_2)/G^0$. Если элемент из $G/G^0$ оставляет на месте некоторую орбиту действия $G^0\cln U\sm U_2$, то найдётся представитель $g\in G$ соответствующего смежного класса, переводящий в себя некоторый вектор этой орбиты: $gv=v$,~$v\in U\sm U_2$.
\end{proof}

Теперь докажем теорему~\ref{Mich}. Пусть группа $G$ конечна. Напомним, что осталось доказать первое утверждение: если $(V\oplus\R^d)/G$\т гладкое многообразие для некоторого $d$, то $G$ порождена псевдоотражениями. Достаточно доказать это для $d=0$: ранг оператора $E-g$ ($g\in G$) на пространствах $V$ и $V\oplus\R^d$ один и тот же.

Итак, предположим, что $V/G$\т гладкое многообразие. Согласно следствию~\ref{1dim}, группа $G$ не содержит отражений. Для всякого вектора $v\in V$ имеем $\ggt v=0$, $V=N_v=N_v^{G_v}\oplus M_v=V^{G_v}\oplus M_v$, а коразмерность страта в $V/G$, соответствующего $v$, равна
$\dim(M_v/G_v)=\dim M_v$, поскольку~$|G_v|<\bes$.

По теореме~\ref{gener} группа $G$ порождена элементами $g$, оставляющими на месте некоторый вектор $v$, которому соответствует страт в $V/G$ коразмерности не более $2$, то есть для которого~$\dim M_v\le2$. Оператор $g\in G_v$ действует тождественно на подпространстве $V^{G_v}$, значит, $\dim\br{(E-g)V}\le\dim M_v\le2$. Таким образом, $g$\т тождественный оператор либо псевдоотражение (отражений в группе $G$ нет).

Тем самым теорема~\ref{Mich} доказана.

\begin{imp}\label{refstab} Пусть стабилизатор $G_v$ некоторого вектора $v\in V$ конечен. Тогда
\begin{nums}{-1}
\item если $(V\oplus\R^d)/G$\т гладкое многообразие для некоторого $d$, то $G_v=\ha{G_v\cap\Om}$, а если $G_v=\ha{G_v\cap\Om}$, то фактор $V/G$
локально в точке $\pi(v)$ является многообразием;
\item если $G_v=\ha{G_v\cap\Om}$, а вектор $v$ лежит на единичной сфере $S\subs V$, то фактор $S/G$ локально в точке $\pi(v)$ является многообразием.
\end{nums}
\end{imp}

\begin{proof} Множество всех элементов группы $G_v$, действующих на $N_v$ псевдоотражением или тождественно, совпадает с $G_v\cap\Om$.

Если $v\in S$, то ортогональное дополнение к подпространству $\ggt v$ в касательном пространстве $T_v(S)=\ha{v}^{\perp}$ есть $N_v(S):=N_v\cap(T_vS)=N_v\cap\ha{v}^{\perp}$. Ясно, что $N_v=N_v(S)\oplus\ha{v}$. Каждый элемент из $G_v\cap\Om$ действует на прямой $\ha{v}$ тождественно, а на $N_v(S)$, как и на $N_v$,\т псевдоотражением либо тождественно.

Осталось воспользоваться теоремами~\ref{slice} и~\ref{Mich}, а также следствием~\ref{slicd}.
\end{proof}

\begin{theorem}\label{stratk} Предположим, что у орбиты общего положения представления $G\cln V$ гомотопическая группа $\pi_k$ ($k>0$) нетривиальна, а в факторе $V/G$ любой страт, отличный от главного, имеет коразмерность более~$k+2$. Тогда $V/G$ не есть гладкое многообразие.
\end{theorem}

\begin{proof} По условию главный страт в $V$ совпадает с $V\sm U_{k+2}$.

Допустим, что $V/G$\т гладкое многообразие. Тогда существует связная $G$\д инвариантная окрестность нуля $U\subs V$, такая что фактор $U/G$ кусочно-диффеоморфен открытому шару размерности~$\dim(V/G)$. Подмножество $U':=U\sm U_{k+2}$ открыто, связно, $G$\д инвариантно и содержится в главном страте пространства~$V$. Подмножество $U'/G$ в $V/G$ связно, а его гомотопические группы $\pi_1\sco\pi_{k+1}$ тривиальны. Отображение $\pi$ даёт локально тривиальное расслоение при ограничении на главный страт, а значит, и при ограничении на~$U'$.

Окрестность нуля $U$ содержит открытый шар $B$ некоторого радиуса $\de>0$ с центром в нуле. Множество $U'\cap B=B\sm U_{k+2}$ непусто,
$G$\д инвариантно, связно и имеет тривиальные гомотопические группы~$\pi_1\sco\pi_{k+1}$. Пусть $o\in U'\cap B$\т некоторый вектор. Он лежит в главном страте, и по условию группа $\pi_k(Go)$ нетривиальна. Расслоение $U'\to U'/G$ имеет слой $Go$. Рассмотрим участок точной гомотопической последовательности
\equ{\label{complex}
\pi_{k+1}(U'/G)\xra{\rh_{k+1}}\pi_k(Go)\xra{\tau_k}\pi_k(U').}
Поскольку группа $\pi_{k+1}(U'/G)$ тривиальна, гомоморфизм $\tau_k$ является вложением. Далее, $Go\subs U'\cap B$, следовательно, для любого элемента $\ga\in\pi_k(Go)$ элементу $\tau_k(\ga)\in\pi_k(U')$ соответствует $k$\д мерный сфероид в $U'$, целиком лежащий в~$U'\cap B$. Этот сфероид стягивается в точку в подмножестве $U'\cap B$ с тривиальной группой $\pi_k$, а значит, и во всём множестве $U'$, откуда элемент $\tau_k(\ga)$ тривиален. Ввиду произвольности $\ga$ гомоморфизм $\tau_k$ тривиален, но в то же время он является вложением, а группа $\pi_k(Go)$ нетривиальна. Получили противоречие.
\end{proof}

\begin{imp}\label{stratd} В условиях теоремы~\ref{stratk} фактор $(V\oplus\R^d)/G$ ни при каком $d$ не является гладким многообразием.
\end{imp}

\begin{proof} Достаточно доказать, что при $d\ge0$ представление~\eqref{act} удовлетворяет условиям теоремы~\ref{stratk}.

Фактор этого действия гомеоморфен $V/G\times\R^d$, страты в нём суть в точности декартовы произведения стратов в $V/G$ на $\R^d$, и их коразмерности равны коразмерностям соответствующих стратов в~$V/G$. Наконец, орбита общего положения для действия $G\cln V$ будет орбитой общего положения и для представления~\eqref{act}.
\end{proof}

\begin{lemma}\label{pitriv} Допустим, что $V$ есть прямая сумма двух $G$\д инвариантных подпространств $V'$ и $V\"$, а фактор действия $G$ на единичной сфере $S\subs V$ является многообразием. Кроме того, предположим, что фактор единичной сферы в $V\"$ по действию $G$\т то есть
$(S\cap V\")/G$\т гомеоморфен замкнутому шару (некоторой размерности), а у фактора $\br{V'\sm\{0\}}/G$  все гомотопические группы тривиальны (под этим также подразумеваем то, что указанный фактор связный). Тогда $V/G$\т многообразие.
\end{lemma}

\begin{proof} Через $\ta_t$ ($t\in\R$, $0\le t\le1$) обозначим линейное отображение в $V$, действующее на $V'$ тождественно, а на $V\"$\т умножением на~$t$. Тогда $\ta\cln[0;1]\times V\to V$\т непрерывное отображение, $\ta_1=E$, $\ta_0(V)=V'$. Из равенств $GV'=V'$ и $GV\"=V\"$ следует, что все
$\ta_t$ коммутируют со всеми операторами из $G$, а множество $U:=V\sm V\"$ инвариантно относительно всех отображений $\ta_t$ и $g\in G$, причём под действием $\ta_0$ оно переходит в~$U\cap V'$. Значит, на самом деле имеются отображения $\wt{\ta}_t\cln U/G\to U/G$, такие что $\wt{\ta}\cln[0;1]\times U/G\to U/G$\т непрерывное отображение, $\wt{\ta}_1$\т тождественное преобразование в $U/G$, а образ $\wt{\ta}_0$ лежит в подмножестве $(U\cap V')/G$, на котором, в свою очередь, все $\wt{\ta}_t$ тождественны. Отсюда фактор $U/G$ гомотопически эквивалентен фактору
$(U\cap V')/G=\br{V'\sm\{0\}}/G$, у которого, по условию, все гомотопические группы тривиальны. Поэтому у $U/G$ также все гомотопические группы тривиальны.

Если из фактора $S/G$ удалить подмножество $(S\cap V\")/G$, гомеоморфное замкнутому шару, то оставшееся подмножество при декартовом умножении на $\R_{>0}$ даёт топологическое пространство, гомеоморфное $U/G$. Значит, у оставшегося топологического пространства все гомотопические группы, как и у $U/G$, тривиальны. Итак, при удалении из многообразия $S/G$ подмножества, гомеоморфного замкнутому шару, остаётся топологическое пространство, у которого все гомотопические группы тривиальны. Следовательно, фактор $S/G$ гомеоморфен сфере, а конус $V/G$ над ним\т векторному пространству.
\end{proof}

На протяжении дальнейшей части работы будем считать группу $G^0$ коммутативной.

В каждой изотипной компоненте $V_{\la}$ ($\la\in P$) представления $G^0\cln V$ стабилизатор любого ненулевого вектора в группе $G^0$ совпадает с ядром гомоморфизма $\la\cln G^0\to\T$. Если вектор $v\in V$ имеет ненулевые проекции на все~$V_{\la}$, то его стабилизатор в группе $G^0$ совпадает с пересечением ядер всех весов $\la\cln G^0\to\T$, которое тривиально; отсюда~$|G_v|<\bes$. Таким образом, типичный вектор пространства $V$ имеет конечный стабилизатор в $G$, типичная орбита имеет размерность $\dim G$, и $\dim(V/G)=\dim V-\dim G$.

Из коммутативности группы $G^0$ следует, что $\Ad(G^0)=\{E\}$. В частности, для любого $v\in V$ группа $G_v\cap G^0$ тождественно действует
на~$\ggt v$.

Для произвольного подмножества $Q\subs P$ рассмотрим группу Ли $\T^{|Q|}$, в которой каждый элемент $z$ определяется своими координатами $z_{\la}\in\T$ ($\la\in Q$). В алгебре же $\Lie\T^{|Q|}=\R^{|Q|}$ каждый элемент $x$ определяется своими координатами $x_\la\in\R$ ($\la\in Q$).

Построим гомоморфизм групп Ли $\ph_Q\cln G^0\to\T^{|Q|}$, $\br{\ph_Q(g)}_{\la}:=\la(g)$. Этот гомоморфизм имеет дифференциал
$d\ph_Q\cln\ggt\to\R^{|Q|},\;\br{d\ph_Q(\xi)}_{\la}=\la(\xi)$.

\begin{prop}\label{1st} Если $V/G$\т многообразие, то множество $P\subs\ggt$ является $1$\д устойчивым.
\end{prop}

\begin{proof} Требуется доказать, что при удалении (с учётом кратности) из $P$ любого веса $\la$ остаётся множество весов с тривиальным пересечением ядер\т что равносильно, с линейной оболочкой~$\ggt$. Если $\la=0$, то доказывать нечего. Если же $\la$\т ненулевой вес, то ему соответствует двумерное неприводимое подпредставление $W\subs V$ группы~$G^0$. В силу следствия~\ref{2dim} группа $G[W]$ конечна, а пересечение ядер всех весов, оставшихся при удалении (с учётом кратности) $\la$ из $P$, есть $\ggt[W]=\Lie G[W]=0$, что и требовалось.
\end{proof}

\begin{lemma}\label{k0} Пусть $W\subs V$\т некоторое $G^0$\д инвариантное подпространство, а группа $G[W]$ имеет размерность~$k>0$. Если при этом $V/G$\т многообразие, то $\dim W\ge2(k+1)$.
\end{lemma}

\begin{proof} Линейное представление $G[W]\cln V$ точное, значит, соответствующее ему множество $P'$ весов $\ggt[W]\to\R$ линейно порождает
$k$\д мерное пространство, двойственное к~$\ggt[W]$. Это множество $1$\д устойчиво, так как состоит в точности из ограничений всех весов
$1$\д устойчивого множества $P$ на подпространство~$\ggt[W]\subs\ggt$. Поскольку $k>0$, в $P'$ содержится по крайней мере $k+1$ ненулевых весов. При этом группа $G[W]$ тождественно действует на подпространстве $W^{\perp}$, отсюда все ненулевые веса из $P'$ соответствуют (с учётом кратностей) двумерным неприводимым подпредставлениям в представлении~$\br{G[W]}^0\cln W$. Поэтому~$\dim W\ge2(k+1)$.
\end{proof}

До конца этого пункта будем предполагать, что $(V\oplus\R^d)/G$\т гладкое многообразие для некоторого~$d$.

\begin{theorem}\label{finstab} Факторгруппа $G/G^0$ порождена смежными классами, содержащими представителя некоторой конечной группы $G_v$,~$v\in V$.
\end{theorem}

\begin{proof} Докажем, что если вектору $v\in V$ соответствует страт в $V/G$ коразмерности не более $2$, то~$|G_v|<\bes$. В самом деле, группа $G_v\cap G^0$ действует тождественно на подпространствах $N_v^{G_v}$ и $\ggt v$, откуда $G_v\cap G^0\subs G[M_v]$. Если $k:=\dim G_v>0$, то
$\dim G[M_v]\ge\dim(G_v\cap G^0)=k$, значит, согласно лемме~\ref{k0}, $\dim M_v\ge2(k+1)$. Поэтому коразмерность соответствующего страта в $V/G$ равна
\equ{\label{codim}
\dim(M_v/G_v)=\dim M_v-\dim G_v\ge2(k+1)-k=k+2>2,}
что противоречит предположению. Следовательно, вектор $v$ имеет конечный стабилизатор.

Теперь можно воспользоваться теоремой~\ref{gener} и получить нужное утверждение.
\end{proof}

\begin{imp}\label{fst1} Факторгруппа $G/G^0$ порождена смежными классами, пересекающимися с~$\Om$. Что равносильно, группа $G$ порождена множеством~$G^0\cup\Om$.
\end{imp}

\begin{proof} Сразу вытекает из теоремы~\ref{finstab} и следствия~\ref{refstab}.
\end{proof}

\begin{imp}\label{fst2} Группа $\Ad(G)$ порождается операторами $\Ad(g)$ по всем~$g\in\Om$.
\end{imp}

\begin{theorem}\label{strat3} Если $\dim G>0$, то в факторе $V/G$ найдётся страт положительной коразмерности, не большей~$3$.
\end{theorem}

\begin{proof} Достаточно доказать, что орбита общего положения для действия $G\cln V$ неодносвязна; после этого, полагая $k:=1$, можно применить
следствие~\ref{stratd}.

Если $H\subs G$\т стабилизатор общего положения, то $|H|<\bes$, а орбита общего положения гомеоморфна многообразию $G/H$ левых смежных классов $G$ по~$H$. Отображение $G\to G/H,\,g\to gH$ является накрытием с конечным слоем $H$, причём накрывающее пространство неодносвязно, так как $G^0$\т тор положительной размерности. Поэтому база $G/H$ тоже неодносвязна.
\end{proof}

\section{Базовые примеры}\label{examp}

\subsection{Конечные группы, порождённые комплексными отражениями}\label{ShSh}

Предположим, что пространство $V$ имеет комплексную структуру, $\dim_{\Cbb}V=n$, и на нём действует компактная группа Ли $G$ линейных над $\R$ операторов, каждый из которых линеен или антилинеен над~$\Cbb$.

\begin{lemma}\label{RCbb} Пусть группа Ли $G$ действует в некотором комплексном пространстве линейными и антилинейными над $\Cbb$ операторами. Тогда для любого $G$\д инвариантного комплексного подпространства найдётся дополнительное $G$\д инвариантное комплексное подпространство.
\end{lemma}

\begin{proof} Из условия следует, что группа $G$ нормализует группу $\Ga:=\{\pm E;\pm I\}$, порождённую оператором $I$ умножения на $i$, поэтому $G\Ga$\т компактная группа Ли, и любое её вещественное линейное представление вполне приводимо. Остаётся заметить, что $(G\Ga)$\д инвариантные вещественные подпространства суть в точности $G$\д инвариантные комплексные.
\end{proof}

Пусть $H\nl G$\т некоторая конечная нормальная подгруппа, порождённая (комплексными) отражениями пространства~$V$.

\begin{theorem}\label{ShShT} Фактор $V/H$ диффеоморфен пространству $V$, причём всякий оператор из $G$, линейный (антилинейный) над $\Cbb$, действует на $V/H$ линейно (антилинейно) над~$\Cbb$.
\end{theorem}

\begin{proof} Группа $G$ действует на линейном пространстве $\Cbb[V]$ по формуле
\eqn{\label{gfv}
gf(v):=\case{
f(g^{-1}v),&g\in\GL_{\Cbb}(V);\\\\
\ol{f(g^{-1}v)},&g\notin\GL_{\Cbb}(V).}}
Линейный (антилинейный) над $\Cbb$ оператор из группы $G$ действует на $\Cbb[V]$ линейно (антилинейно) над~$\Cbb$. Кроме того, $g(f_1f_2)=(gf_1)(gf_2)$ при $g\in G$,~$f_i\in\Cbb[V]$.

Первое утверждение теоремы следует из теоремы Шепарда-Тодда-Шевалле: отображение факторизации можно задать формулой
$\pi_H\cln V\to V/H,\,v\to\br{f_1(v)\sco f_n(v)}$, где $\{f_1\sco f_n\}$\т любая система алгебраически независимых однородных образующих алгебры
инвариантов $B:=\Cbb[V]^H$. Поскольку $H\nl G$, имеем~$GB=B$. Пусть $B_k\subs B$\т подпространство однородных инвариантов степени~$k$.

В каждом $B_k$ ($k>0$) можно выделить подпространство $C_k$, натянутое на все $B_iB_{k-i}$ ($0<i<k$), то есть порождённое всевозможными произведениями многочленов из $B$ степеней менее $k$, которые (произведения) сами однородны и имеют степень~$k$. Конечномерное комплексное пространство $B_k$ и его комплексное подпространство $C_k$ инвариантны относительно $G$, значит, существует $G$\д инвариантное комплексное подпространство $C'_k\subs B_k$, такое что $B_k=C_k\oplus C'_k$. Выбрав в каждом $C'_k$ базис, мы и получим систему $(f_1\sco f_n)$ алгебраически независимых однородных образующих алгебры~$B$. Зададим отображение факторизации формулой $\pi_H(v):=\br{f_1(v)\sco f_n(v)}$. При этом $\pi_H(0)=0$, так как $f_i$\т многочлены без свободного члена.

Для любого $g\in G$ имеем $gC'_k=C'_k$. Если при этом $g$\т линейный (антилинейный) оператор в $V$, то он действует на $V/H$ линейно (антилинейно) над $\Cbb$, поскольку каждое число $f_i(gv)$ (число $\ol{f_i(gv)}$) комплексно-линейно выражается через числа $f_j(v)$, такие что~$\deg f_j=\deg f_i$.
\end{proof}

\begin{imp} Для любого $d\ge0$ фактор $(V\oplus\R^d)/H$ диффеоморфен пространству $V\oplus\R^d$, причём всякий линейный (антилинейный) над $\Cbb$ оператор из $G$ действует на прямом слагаемом $\R^d$ этого фактора тождественно, а на прямом слагаемом $V$\т линейно (антилинейно) над~$\Cbb$.
\end{imp}

Как видно из доказательства теоремы~\ref{ShShT}, в пространстве $V/H$ линейные оболочки всех базисных векторов, соответствующих однородным образующим какой-либо фиксированной степени, $G$\д инвариантны. В частности, если представление $G\cln(V/H)$ неприводимо, то степени всех однородных образующих равны.

Итак, мы получили линейное (над $\R$) представление группы Ли $G/H$ в пространстве~$V/H$. Докажем, что оно точно, то есть что любой элемент $g\in G$, действующий на $V/H$ тождественно, лежит в~$H$. В самом деле, для любого $v\in V$ имеем $gv\in Hv$, и $h^{-1}gv=v$ для некоторого~$h$. Мы видим, что $V$ есть объединение линейных (над полем $\R$) подпространств неподвижных точек линейных (тоже над $\R$) операторов $h^{-1}g$,~$h\in H$. Так как $|H|<\bes$, для некоторого $h\in H$ это подпространство совпадает с $V$, откуда~$g=h$. Группа $G/H$ может быть отождествлена со своим образом
в $\GL_{\R}(V/H)$ при полученном точном представлении.

Пусть в группе $G/H$ есть некоторая конечная порождённая отражениями подгруппа $G'/H$, где $H\subs G'\subs G$. В частности, $G'/H$ состоит только из линейных над $\Cbb$ операторов в~$V/H$. Значит, все элементы из $G'$ действуют на $V$ комплексно-линейно. В силу конечности групп $H$ и $G'/H$, группа $G'$ также конечна. Поскольку $G'/H$ порождена отражениями, $V/G'\cong(V/H)/(G'/H)\cong V/H\cong V$, отсюда группа $G'$ порождена отражениями пространства~$V$. В частности, если $H$ порождена \textit{всеми} отражениями группы $G$, то $G/H$ не содержит отражений, иначе в $G/H$ была бы нетривиальная конечная порождённая отражениями подгруппа $G'/H$, $H\subsneq G'\subs G$, и группа $G'$ была бы порождена отражениями, что невозможно.

Допустим, что $V_0=0$, а подгруппа $G^0\subs G$ центральна, то есть $\Ad(G)=\{E\}$. Тогда $G$ оставляет на месте все изотипные компоненты $V_{\la}$ ($\la\in P$) представления $G^0\cln V$ и действует на каждой из них комплексно-линейно. Таким образом, на пространстве $V$ как на прямой сумме $V_{\la}$ тоже вводится комплексная структура, причём группа $G$ действует линейно на комплексном пространстве~$V$. Для одномерной группы Ли $G$ условие $V_0=0$ равносильно тому, что стабилизатор любого ненулевого вектора в $V$ конечен (что то же самое, не содержит~$G^0$).

\begin{prop}\label{fin} Если $\dim G=1$, $V_0=0$, $\Ad(G)=\{E\}$, $\dim_{\Cbb}V>1$, то подгруппа $H\nl G$, порождённая всеми (комплексными) отражениями группы $G$, конечна.
\end{prop}

\begin{proof} Так как $GV_{\la}=V_{\la}$ для всех $\la$, любое отражение в $G$ действует тождественно на всех $V_{\la'}$, кроме одного~$V_{\la}$, на котором\т комплексным отражением, то есть лежит в группе~$G[V_{\la}]$. Для каждого $\la$ обозначим через $H_{\la}$ группу, порождённую всеми отражениями из~$G[V_{\la}]$. Тогда $H$ есть прямое произведение всех~$H_{\la}$.

При наличии по крайней мере двух изотипных компонент каждая $H_{\la}$ тождественно действует на всех $V_{\la'}$, $\la'\ne\pm\la$, поэтому
$H_{\la}\cap G^0$ содержится в ядре нетривиального веса $\la'\cln G^0\to\T$. Это ядро конечно, отсюда и группа $H_{\la}$ конечна. Значит, конечна и группа $H$ как прямое произведение всех~$H_{\la}$.

Если же есть только одна изотипная компонента $V_{\la}=V$, то $G^0$ действует на пространстве $V$ скалярными операторами, а для произвольного отражения $g\in G$ оператор $g^{|G/G^0|}\in G^0$\т отражение либо тождественный. Поскольку в пространстве размерности более $1$ скалярный оператор не может быть отражением, $g^{|G/G^0|}=E$,
\eqn{\label{det}
(\det g)^{|G/G^0|}=1}
(определитель комплексный). Все элементы $g\in G$, удовлетворяющие~\eqref{det}, образуют подгруппу, в которой содержатся все отражения из $G$, а значит, и подгруппа~$H$. Но любой оператор $g\in H\cap G^0$ является скалярным, имеет вид $\la E$, причём, ввиду~\eqref{det}, число $\la$ может принимать лишь конечное количество значений. Следовательно, $|H\cap G^0|<\bes$, $|H|<\bes$.
\end{proof}

Фактор $(V\oplus\R^d)/H$ ($d\ge0$) диффеоморфен пространству $V\oplus\R^d$, одномерная группа $G/H$ действует на его прямом слагаемом $V$ комплексно-линейными операторами, среди которых нет отражений, а на прямом слагаемом $\R^d$\т тождественно. Связная компонента единицы группы $G$ лежит в её центре, значит, этим же свойством обладает и её гомоморфный образ~$G/H$. Наконец, для действия $(G/H)\cln(V/H)$ стабилизатор любого ненулевого вектора конечен. В самом деле, иначе группа $(G/H)^0=(G^0H)/H$ оставляла бы на месте ненулевой вектор $v'=\pi_H(v)\in V/H$, где $v\in V$,
$v\ne0$, тогда $G^0v\subs Hv$, и, в силу связности $G^0$ и конечности $H$,~$G^0v=\{v\}$. Последнее противоречит условию~$V_0=0$.

\begin{theorem}\label{ab} Если $\dim G=1$, $\Ad(G)=\{E\}$, $V_0=0$, а множество $P$ содержит не менее трёх ненулевых весов, то ни при каком $d$ фактор $(V\oplus\R^d)/G$ не является гладким многообразием.
\end{theorem}

\begin{proof} Предположим, что $(V\oplus\R^d)/G$\т гладкое многообразие, где~$d\ge0$.

Группа $G$ линейно действует на $V$ как на комплексном пространстве. По условию $\dim_{\Cbb}V>2$. Поэтому все элементы из $G$, действующие на $V$ комплексными отражениями, порождают конечную подгруппу~$H\nl G$. Как уже известно, фактор $(V\oplus\R^d)/H$ диффеоморфен пространству~$V\oplus\R^d$. На последнем канонически действует группа $G/H$, и фактор этого действия является гладким многообразием, так как $(V\oplus\R^d)/G$\т гладкое многообразие.

Тем самым всё сведено к случаю, когда $(V\oplus\R^d)/G$\т гладкое многообразие, группа $G$ действует на $V$ комплексно-линейно без комплексных отражений, причём по-прежнему $V_0=0$,~$\dim_{\Cbb}V>2$. По теореме~\ref{strat3} в $V/G$ есть страт, отличный от главного, коразмерности не более~$3$. Рассмотрим в таком страте точку $\pi(v)\in V/G$, где~$v\in V$. Если $|G_v|=\bes$, то $v=0$, поэтому указанный страт в $V/G$ состоит из одной точки $\pi(0)$ и имеет коразмерность $\dim(V/G)=2\dim_{\Cbb}V-1\ge5$. Значит, $|G_v|<\bes$, и коразмерность страта в $V/G$ равна $\dim M_v$, откуда
$0<\dim M_v<4$. Поскольку $\Ad(G)=\{E\}$, комплексное подпространство $V^{G_v}$ совпадает с $N_v^{G_v}\oplus\ggt v$ и имеет комплексную коразмерность $\frac12\dim M_v<\frac42=2$. В то же время группа $G_v$ нетривиальна, так как соответствующий страт не главный, и не содержит комплексных отражений. Получили противоречие.
\end{proof}

\subsection{Действие тора на комплексном пространстве}\label{torum}

Допустим, что пространство $V$ вновь имеет комплексную структуру, $\dim_{\Cbb}V=n\in\N$, и разлагается в прямую сумму попарно ортогональных одномерных комплексных подпространств $W_1\sco W_n$.

\begin{theorem}\label{stor} Пусть в пространстве $V$ действует $n$\д мерная компактная группа Ли $G\subs\Or(V)$, переставляющая подпространства $W_i$, причём если $n>2$, то $GW_i=W_i$ для всех~$i$. Тогда фактор единичной сферы $S\subs V$ по действию группы $G$ гомеоморфен замкнутому $(n-1)$\д мерному шару; в частности, все гомотопические группы факторов $S/G$ и $\br{V\sm\{0\}}/G$ тривиальны.
\end{theorem}

\begin{proof} Группа Ли $G^0$ переводит в себя все подпространства $W_i$ и действует на каждом из них комплексно-линейно. Поэтому она вкладывается в
связную группу Ли $\T^n$, значит, это вложение есть изоморфизм, откуда $V/G^0\cong\R_{\ge0}^n$. Пусть $g\in G$\т произвольный элемент. Если он переводит в себя все подпространства $W_i$, то действует на факторе $V/G^0$ тождественно, поскольку принадлежит~$\Or(V)$. В противном случае имеем $n=2$, $g$ переводит $W_1$ и $W_2$ друг в друга и действует на факторе $V/G^0\cong\R_{\ge0}^2$ перестановкой двух координат по~$\R_{\ge0}$.
В любом случае фактор $V/G$ гомеоморфен замкнутому $n$\д мерному полупространству, причём $\pi(0)$\т его граничная точка, а фактор $S/G$ гомеоморфен замкнутому $(n-1)$\д мерному шару: либо это симплекс на единичной сфере в $\R^n$, высекаемый координатными гиперплоскостями, либо некоторая дуга (одномерной) окружности.
\end{proof}

Предположим, что $n>1$. Выберем в каждом пространстве $W_i$ вектор $e_i$ единичной длины. Пусть в пространстве $V=\ha{e_1\sco e_n}_{\Cbb}$ действует $(n-1)$\д мерная группа Ли $H\subs\Or(V)$ всех комплексно-линейных диагональных операторов $\diag(\la_1\sco\la_n)$, таких что $|\la_i|=1$ и
$\la_1^{a_1}\sd\la_n^{a_n}=1$, $a_i\in\N$. Группу $H$ можно отождествить с подтором в торе $\T^n$, задаваемым одним характером с натуральными показателями степеней $a_1\sco a_n$.

\begin{theorem}\label{tor} Фактор $V/H$ диффеоморфен вещественному векторному пространству $\R^{n+1}$, причём всякий ортогональный оператор
$g\in N(H)$ действует на нём линейно.
\end{theorem}

\begin{proof} Пространства $W_i$ суть одномерные комплексные представления группы $H$. Если два из них: $W_i$ и $W_j$\т изоморфны как её вещественные представления, то соответствующие им веса (гомоморфизмы $H\to\T$) совпадают с точностью до поточечного обращения. Это означает, что характер $(\la_1\sco\la_n)\to\la_i\la_j^{\pm1}$ тора $\T^n$ тривиален на подторе $H$, то есть его показатели степеней имеют вид $ca_1\sco ca_n$, где $c\in\Z$, откуда $n=2$ и $a_1=a_2=1$.

Вначале предположим, что вещественные представления $W_i$ группы $H$ попарно неизоморфны. Тогда они являются изотипными компонентами представления $H\cln V$ и потому переставляются любым оператором из $N(H)$.

Выясним, как устроен фактор $V/H$. Если два вектора пространства $V$ лежат в одной орбите действия $H\cln V$, то, очевидно, модули их координат попарно равны, и произведения их координат в степенях $a_i$ также равны. Если имеются два вектора с попарно равными модулями координат и среди этих модулей есть нулевой, то они лежат в одной орбите. Действительно, диагональный элемент, соответствующий нулевой координате, может быть выбран произвольно, поэтому вначале выберем все, кроме него (здесь тоже возможна неоднозначность), а затем\т оставшийся, руководствуясь условием
$\la_1^{a_1}\sd\la_n^{a_n}=1$. Если модули координат попарно равны и все отличны от нуля, то два вектора лежат в одной орбите тогда и только тогда, когда произведения координат в степенях $a_i$ равны. Отображение факторизации $\pi_H\cln V\to\R^{n+1}=\R^{n-1}\oplus\Cbb$ построим так: всякому вектору $z=(z_1\sco z_n)$ сопоставим вектор $\br{|z_1|^2\sco|z_n|^2}\in\R^n$, который спроецируем на гиперплоскость $\sumiun x_i=0$ вдоль вектора $(1\sco1)$, то есть ортогонально. Полученная проекция будет, по определению, координатой вектора $\pi_H(z)$ в $\R^{n-1}$; координатой же в $\Cbb$ будет число $z_1^{a_1}\sd z_n^{a_n}$. Построенное отображение $\pi_H$, очевидно, гладко и постоянно на орбитах $H$. С другой стороны, если известна
координата в $\R^{n-1}$ (а именно, набор квадратов модулей координат с точностью до прибавления одного и того же числа), а также произведение этих модулей в степенях $a_i$, то сами модули восстанавливаются однозначно: при строгом увеличении всех модулей их произведение в степенях $a_i$ также строго увеличивается, крайние значения\т $0$ и $+\bes$. Отсюда следует, что отображение $\pi_H$ сюръективно и разделяет орбиты.

Теперь рассмотрим произвольный элемент $g\in N(H)$. Он, как мы уже знаем, переставляет подпространства~$W_i$. Его можно представить в виде композиции $g=g_1\circ g_2\circ g_3$ линейных над $\R$ преобразований: $g_1$ осуществляет некоторую перестановку $\si$ комплексных координат, то есть $(g_1z)_i:=z_{\si(i)}$; $g_2$ сопрягает некоторые координаты, сохраняя остальные; $g_3$ умножает координату $z_i$ на число $\nu_i$, $|\nu_i|=1$. Очевидно, что $g_3Hg_3^{-1}=H=gHg^{-1}$, поэтому оператор $g_0:=g_1g_2$ тоже лежит в нормализаторе $H$. Но $g_0Hg_0^{-1}$\т это группа всех
диагональных операторов $\diag(\la_1\sco\la_n)$, таких что $|\la_i|=1$ и $\la_1^{b_1}\sd\la_n^{b_n}=1$, где $b_i=a_{\si(i)}$ ($b_i=-a_{\si(i)}$), если координата с номером $\si(i)$ сохраняется (сопрягается) оператором $g_2$. А так как $g_0Hg_0^{-1}=H$, числа $b_1\sco b_n$ пропорциональны числам $a_1\sco a_n$, откуда
\eqn{\label{asi}
a_i=a_{\si(i)}\,\fa i,}
а оператор $g_2$ либо тождественный, либо сопрягает все координаты.

Посмотрим теперь, как такой элемент $g$ будет действовать на факторе. Координата в $\R^{n-1}$ меняется лишь под действием $g_1$, причём в этом случае модули координат переставляются, в $\R^n$ происходит линейное преобразование перестановки координат, и его надо лишь ограничить на инвариантную гиперплоскость $x_1\spl x_n=0$. Что же касается координаты в $\Cbb$, то она умножается на $\nu_1^{a_1}\sd\nu_n^{a_n}$ при $g_3$, сохраняется либо сопрягается при~$g_2$, а при $g_1$ сохраняется, согласно~\eqref{asi}. Итак, элемент $g$ действует на $V/H$ линейно.

Теперь предположим, что среди вещественных представлений $H\cln W_i$ по крайней мере два изоморфны. В этом случае, как мы знаем, $n=2$ и $a_1=a_2=1$.

Если обозначить вектор $e_1$ символом $\idb$, $ie_1$\т символом $\ib$, вектор $e_2$\т $\jb$, а $ie_2$\т $\kb$, то пространство $V$ отождествится с пространством $\Hbb$ кватернионов, причём левое умножение на $\ib$ здесь, как легко видеть, будет совпадать с умножением на $i$ в пространстве над $\Cbb$. Если $v\in\ha{\jb,\kb}$, то $v\ib=-\ib v$, откуда $vz=\ol{z}v$ при $z\in\Cbb$ (комплексные числа отождествляем с кватернионами из
$\ha{\idb,\ib}$). Таким образом, левое умножение на комплексное число $z$\т это в точности умножение на элемент $z$ поля $\Cbb$ в пространстве над этим полем, а правое умножение на $z\in\Cbb$\т это диагональный оператор $\diag(z,\ol{z})$ в комплексном пространстве. Следовательно, группа $H$ действует правыми умножениями на комплексные числа с модулем $1$. Рассмотрим отображение $\pi_H\cln\Hbb\to\Hbb,\,v\to v\ib\ol{v}$ (сопряжение здесь
кватернионное\т замена знака всех мнимых частей и сохранение вещественной). Легко понять, что $\pi_H(\Hbb)=\ha{\ib,\jb,\kb}$, а $\pi_H(v_1)=\pi_H(v_2)$ тогда и только тогда, когда $v_1=v_2z$, $z$\т комплексное число с модулем $1$. Итак, слои отображения
$\pi_H$ суть орбиты группы $H$, значит, фактор $V/H$ диффеоморфен $\R^3$.

Нормализатор подгруппы $H$ в группе ортогональных операторов порождается $H$, левыми умножениями на кватернионы с нормой $1$ и правым умножением на $\jb$; на факторе $V/H$ он действует как группа $\Or_3(\R)$.
\end{proof}

Если $n=2$, то $\dim H=1$, а для любого ортогонального оператора $g\in N(H)$ оператор $\Ad_{\hgt}(g)$ на одномерной алгебре $\hgt:=\Lie H$ равен
$(\pm E)$.

\begin{stm} Пусть $n=2$, а группа $H$ связна. Тогда ортогональный оператор $g\in N(H)$ действует на $V/H\cong\R^3$ отражением, если и только если $g^2=E$ и $\Ad_{\hgt}(g)=-E$.
\end{stm}

\begin{proof} Если $\Ad_{\hgt}(g)=\ep E$, $\ep=\pm1$, то определитель оператора $g$ на пространстве $V/H$ равен $\ep$. Поэтому можно считать, что $\Ad_{\hgt}(g)=-E$. В таком случае в пространстве $V$ любой оператор из $\hgt\sm\{0\}$ невырожден и антикоммутирует с $g$, откуда
$\dim\br{\Ker(E-g)}=\dim\br{\Ker(E+g)}=\frac{1}{2}\dim\Br{\Ker\br{E-g^2}}$. Кроме того, всякая орбита группы $H$, инвариантная относительно $g$,
имеет непустое конечное пересечение с подпространством $\Ker(E-g)$. Значит, размерности подпространств неподвижных точек для действия $g$ на пространствах $V$ и $V/H$ совпадают, то есть $\dim(V/H)^g=\dim\br{\Ker(E-g)}=\frac{1}{2}\dim\Br{\Ker\br{E-g^2}}$. В частности,
$\dim(V/H)^g=2\Lra g^2=E$.
\end{proof}

\begin{lemma}\label{O3fin} Пусть $G$\т конечная подгруппа в $\Or_3(\R)$, содержащая хотя бы одно отражение. Тогда для её тавтологического представления в $V:=\R^3$ все гомотопические группы фактора $\br{V\sm\{0\}}/G$ тривиальны.
\end{lemma}

\begin{proof} Достаточно доказать, что все гомотопические группы фактора единичной сферы $S\subs V$ по действию $G$ тривиальны.

Из условия следует, что подгруппа $H\subs G$, порождённая всеми отражениями группы $G$, нетривиальна. Поэтому для её действия на $S$ существует
фундаментальная область $M\ne S$\т замкнутый выпуклый сферический многоугольник, гомеоморфный кругу. Имеем $G=H\lhdp\De$, где
$\De:=\{g\in G\cln gM=M\}$. 
Группа $\De$ не содержит отражений и оставляет на месте центр сферического многоугольника $M$, а значит, действует на $M$ вращениями. Фактор $M/\De$ также гомеоморфен кругу, в частности, все его гомотопические группы тривиальны.

Никакие две точки из $M$ не лежат в одной орбите действия $H\cln S$, откуда фактор $S/G$ гомеоморфен $M/\De$, и все его гомотопические группы тривиальны.
\end{proof}

\begin{theorem}\label{dtor} Пусть одномерная группа $G$ действует на четырёхмерном пространстве $V$, а множество $P$ не содержит нулевых весов. Также предположим, что найдётся инволюция $g\in G$, такая что $\Ad(g)=-E$. Тогда все гомотопические группы фактора $\br{V\sm\{0\}}/G$ тривиальны.
\end{theorem}

\begin{proof} Поскольку $\dim G=1$, группа $\Ad(G)$ совпадает с $\{\pm E\}$.

Пространство $V$ можно разложить в прямую сумму двух неприводимых подпредставлений группы $G^0$, каждое из которых имеет структуру одномерного комплексного пространства. Значит, группа $G^0$ вкладывается в двумерный тор $\T^2=\T\times\T$, причём образ при вложении не содержит ни один из прямых сомножителей и потому задаётся одним характером с ненулевыми (можно считать, что натуральными) показателями степеней. Следовательно, $V$ есть
двумерное комплексное пространство, на котором $G^0$ действует всеми диагональными комплексно-линейными операторами $\diag(\la_1,\la_2)$, такими что $|\la_i|=1$ и $\la_1^{a_1}\la_2^{a_2}=1$, $a_i\in\N$. По теореме~\ref{tor} фактор $V/G^0$ гомеоморфен трёхмерному пространству $\R^3$, отображение факторизации переводит нулевой вектор в нулевой, а $G$ действует на $V/G^0$ линейно. Инволюция $g\in G$, для которой $\Ad(g)=-E$, действует на $V/G^0$ отражением. Имеем $V/G\cong(V/G^0)/(G/G^0)$, $|G/G^0|<\bes$, и осталось применить лемму~\ref{O3fin}.
\end{proof}

\section{Системы векторов конечномерного пространства}\label{vect}

Изучим более подробно устойчивость конечных множеств векторов в пространстве. Пусть в некотором конечномерном векторном пространстве над полем $\F$ имеется конечное множество $P$ векторов. Некоторые из этих векторов могут совпадать; кратность будет учитываться. Количество ненулевых векторов множества $P$ будем обозначать через $\hn{P}$.

\begin{df}\label{triv} Множество векторов будем называть \textit{тривиальным}, если все его векторы нулевые.
\end{df}

При $q>0$ всякое $q$\д устойчивое множество $(q-1)$\д устойчиво. В свою очередь, $(q-1)$\д устойчивое множество $P$ является $q$\д устойчивым тогда и только тогда, когда при удалении из него любых $q$ векторов остаётся множество с линейной оболочкой $\ha{P}$.

Как было установлено в п.~\ref{introd}, если множество $P$ нетривиально и $q$\д устойчиво, то $\hn{P}\ge\dim\ha{P}+q$. Кроме того, при удалении из $P$ любых $q$ ненулевых векторов линейная оболочка не должна измениться. Следовательно, в случае равенства $\hn{P}=\dim\ha{P}+q$ любые $\dim\ha{P}$ ненулевых векторов линейно порождают пространство $\ha{P}$ и потому линейно независимы; значит, любые ненулевые векторы множества $P$ в количестве не более $\dim\ha{P}$ линейно независимы. Отсюда непосредственно вытекает следующее утверждение.

\begin{stm}\label{indep} В $q$\д устойчивом множестве $P$, удовлетворяющем неравенству $\hn{P}\le\dim\ha{P}+q$, любые ненулевые векторы в количестве не более $\dim\ha{P}$ линейно независимы.
\end{stm}

Если $q$\д устойчивое множество $P$ допускает собственное разложение на две компоненты $P_1$ и $P_2$, то множества $P_i$ нетривиальны и
$q$\д устойчивы. Справедливы неравенства $\hn{P_i}\ge\dim\ha{P_i}+q$, поэтому при $q>0$ имеем
\equ{\label{P12}
\hn{P}=\hn{P_1}+\hn{P_2}\ge\hr{\dim\ha{P_1}+q}+\hr{\dim\ha{P_2}+q}=\dim\ha{P}+2q>\dim\ha{P}+q.}
Из этого можно вывести следующее утверждение.

\begin{stm}\label{irred} Если $q>0$, то любое $q$\д устойчивое множество $P$, удовлетворяющее неравенству $\hn{P}\le\dim\ha{P}+q$, неразложимо.
\end{stm}

Нас будут интересовать в первую очередь $1$\д устойчивые и $2$\д устойчивые множества.

Множество $P$ является $1$\д устойчивым тогда и только тогда, когда любой вектор из $P$ линейно выражается через остальные. Всякое $2$\д устойчивое множество $1$\д устойчиво, но обратное неверно. Приведём один критерий $2$\д устойчивости некоторого $1$\д устойчивого множества векторов $P$.

Рассмотрим произвольные два вектора в множестве $P$ (не обязательно $1$\д устойчивом). Объявим их эквивалентными, если и только если для любого линейного соотношения между векторами из $P$ коэффициенты при этих двух векторах равны или не равны нулю одновременно. Ясно, что полученное отношение между векторами из $P$ является отношением эквивалентности. Оно не меняется при умножении некоторых векторов на ненулевой скаляр. Каждый нулевой вектор в $P$ может быть эквивалентен только самому себе. Поэтому если в классе эквивалентности содержится более одного вектора, то все они
ненулевые. Как легко понять, удаление и добавление нулевых векторов не влияет на эквивалентность в множестве $P$.

\begin{theorem}\label{crit} Произвольное $1$\д устойчивое множество является $2$\д устойчивым тогда и только тогда, когда все его элементы попарно неэквивалентны.
\end{theorem}

\begin{proof}
Предположим, что $P$\т $2$\д устойчивое множество. Покажем, что никакие два вектора $a,b\in P$ не эквивалентны. Действительно, при удалении векторов $a$ и $b$ из $P$ (с учётом кратностей всех векторов) должно оставаться множество $P'$, для которого $\ha{P'}=\ha{P}$. Так как $b\in\ha{P'}$, существует линейное соотношение между векторами из $P$ с ненулевым коэффициентом при $b$ и нулевым\т при $a$. Отсюда следует, что $a$ и $b$ неэквивалентны.

Обратно, пусть все элементы $1$\д устойчивого множества $P$ попарно неэквивалентны. Посмотрим, что произойдёт при удалении из него двух
произвольных векторов $a,b\in P$. Поскольку $a$ и $b$ неэквивалентны, существует линейное соотношение между векторами из $P$, у которого, не умаляя общности, коэффициент при $a$ равен нулю, а при $b$\т не равен (возможно, и наоборот). Тогда при удалении из $P$ элементов $a$ и $b$ остаётся множество $P'$ с линейной оболочкой, содержащей $b$. Значит, $\ha{P'}=\ba{P'\cup\{b\}}=\ba{P\sm\{a\}}$, а так как $P$\т $1$\д устойчивое множество, $\ba{P\sm\{a\}}=\ha{P}$, откуда $\ha{P'}=\ha{P}$.
\end{proof}

Всякий класс эквивалентности $N$ произвольного множества $P$, состоящий более чем из одного вектора, вместе с любым своим вектором $v$ содержит и все его копии из $P$. В самом деле, если $w\in P$\т другая копия вектора $v$, то, ввиду линейного соотношения $v-w=0$, ни один вектор из $P$, кроме $v$ и $w$, не эквивалентен $v$, то есть $N\subs\{v;w\}$. Поскольку $\hn{N}>1$, имеем $w\in N$.

Пусть $P_1$ и $P_2$\т конечные множества векторов в конечномерных пространствах над полем $\F$, причём $P_2$ получается из $P_1$ линейным изоморфизмом $A$ с точностью до ненулевых скаляров: $P_2=\{\ep_vAv\cln v\in P_1\}$, $\ep_v\in \F^*$. Тогда любой класс эквивалентности $N_1$ множества $P_1$ переходит под действием $A$ в некоторый класс эквивалентности $N_2$ множества $P_2$ с точностью до ненулевых скаляров: $\ep_vAv\in N_2\Lra v\in N_1$.

Рассмотрим множество $P$ и класс эквивалентности $N$ в нём. Если $P$\т $1$\д устойчивое множество, то существует линейное соотношение
$\sums{v\in P}c_vv=0$ с ненулевым коэффициентом $c_w$ при некотором векторе $w\in N$. Тогда $\sums{v\in N}c_vv\in\ha{P\sm N}$\т нетривиальная линейная комбинация векторов из $N$, лежащая в линейной оболочке остальных. В любой такой нетривиальной линейной комбинации все коэффициенты $c_v$, $v\in N$, отличны от нуля, поскольку $N$\т класс эквивалентности. Значит, такая нетривиальная линейная комбинация ровно одна с точностью до пропорциональности, с коэффициентами $c_v\ne0$.

В дальнейшем нам будет интересен только случай $\F=\R$, хотя до сих пор это не было нигде использовано. Пусть $N\subs P$\т некоторый класс эквивалентности. Сменив, если нужно, знак у некоторых элементов этого класса, можно добиться того, чтобы нетривиальная линейная комбинация векторов из $N$, лежащая в линейной оболочке остальных, имела коэффициенты $c_v$ ($v\in N$) одного и того же знака\т можно считать, что положительные. Эту процедуру можно проделать независимо с каждым классом $N$, и тогда для всех классов будем иметь $c_v>0$, $v\in N$. При этом множество $P$ не обязано
перейти в себя (знаки сменятся). При смене знаков сохранится $1$\д устойчивость множества, а также порядки всех классов эквивалентности.

\section{Представление с $2$\д устойчивым множеством весов}\label{2stab}

В этом пункте будем предполагать, что для данного действия $G\cln V$ множество $P$ весов, соответствующих неприводимым компонентам с учётом кратностей, $2$\д устойчиво. Как мы уже знаем, на пространстве $\ggt$ существует инвариантное относительно $\Ad(G)$ скалярное умножение.

Если $v\in V$\т произвольный вектор, то множество весов, соответствующих представлению $G_v\cln V$, состоит в точности из ограничений весов $2$\д устойчивого множества $P$ на подпространство $\ggt_v\subs\ggt$ и потому $2$\д устойчиво. Представлению $G_v\cln N_v$ также соответствует
$2$\д устойчивое множество весов, так как группа $G_v\cap G^0$ тождественно действует на $\ggt v$.

Рассмотрим некоторый элемент $g\in\Om$ и положим $A:=\Ad(g)$\т линейный оператор на пространстве $\ggt$, $P^A:=P\cap\Ker(E-A)$, $r:=\rk(E-A)$. Оператор $g$ в пространстве $V$ переводит в себя все изотипные компоненты представления $G^0\cln V$ тогда и только тогда, когда $P\subs\Ker(E-A)\cup\Ker(E+A)$\т что равносильно,
\eqn{\label{Ppm}
P\sm P^A\subs\Ker(E+A).}

\begin{theorem}\label{nocent} Допустим, что $A\ne E$. Тогда
\begin{nums}{-1}
\item множество $(E-A)P$ содержит ровно $r+2$ ненулевых векторов, является неразложимым, а любые его ненулевые векторы в количестве не более $r$
линейно независимы;
\item оператор $g$ действует тождественно на всех изотипных компонентах $V_{\la}$, $\la\in P^A$, в том числе на $V_0$;
\item если множество $P\sm P^A$ не лежит в подпространстве $\Ker(E+A)$, то оно состоит ровно из трёх весов $\la^+,\la^-,\la$, каждый из которых
входит в $P$ с кратностью $1$, причём $A\la=-\la$,
\eqn{\label{Alapm}
A\la^{\pm}=\la^{\mp}\ne\pm\la^{\pm},}
а оператор $A$ есть отражение относительно гиперплоскости $\ha{\la}^{\perp}$.
\end{nums}
\end{theorem}

\begin{proof} Поскольку $A\ne E$, число $r$ положительно. Множество $(E-A)P$, как и $P$, является $2$\д устойчивым. Его линейная оболочка $(E-A)\ha{P}=(E-A)\ggt$ имеет размерность $r>0$, откуда $\bn{(E-A)P}\ge r+2$. Имеем
\eqn{\label{rkPg}
\rk(E-g)=\rk(E-A)+\om(g)\le r+2\le\bn{(E-A)P}.}

Оператор $g$ в пространстве $V$ переставляет изотипные компоненты представления $G^0\cln V$. Рассмотрим один из независимых циклов этой перестановки: $V_{\la},gV_{\la}=V_{A\la}\sco g^kV_{\la}=V_{A^k\la}$, $k\ge0$, $k+1$\т длина цикла, $g^{k+1}V_{\la}=V_{\la}$. Каждое $g^iV_{\la}$ есть прямая сумма
одного и того же количества неприводимых компонент, которое обозначим через $l$ (для разных циклов $l$ может быть разным). Среди весов
$\la,A\la\sco A^k\la\in P$ никакие два не равны и не противоположны, каждый имеет кратность $l$, $A^{k+1}\la=\pm\la$, то есть $A^{k+1}\la=(-1)^{\ep}\la$, где $\ep=0,1$ (если $\la=0$, то положим $\ep:=0$). Пусть $Q\subs P$\т подмножество, состоящее из весов
$\la,A\la\sco A^k\la$, каждый из которых берётся с кратностью $l$. Ясно, что $V_Q=\oplusl{i=0}{k}g^iV_{\la}$, и $gV_Q=V_Q$.

Предположим, что $\la\ne0$. Тогда ядро ограничения на $V_Q$ оператора $E-g$ имеет размерность не более $2l$. В самом деле, если
$(E-g)(v_0\spl v_k)=0$, $v_i\in g^iV_{\la}$, то $gv_i=v_{i+1}$ ($i<k$), $gv_k=v_0$, а также $g^{k+1}v_0=v_0$. Следовательно, если $v_0=0$, то и все $v_i$ равны нулю, то есть ядро проецируется на $V_{\la}$ инъективно. Кроме того, его размерность равна размерности ядра ограничения $E-g^{k+1}$ на $V_{\la}$. Если при этом $\ep=1$, $A^{k+1}\la=-\la$, то ограничение $g^{k+1}$ на $V_{\la}$ антилинейно, поэтому ядра ограничений на $V_{\la}$ операторов $E-g^{k+1}$ и $E+g^{k+1}$ получаются друг из друга умножением на $i$, имеют одинаковую размерность и тривиальное пересечение. Значит,
размерность каждого из них не превосходит $l$. Итак, размерность ядра ограничения $E-g$ на подпространство $V_Q$ не превосходит $2l$, а при $\ep=1$\т не превосходит $l$; в любом случае эта размерность не превосходит $2l-\ep l$. Отсюда
$\dim\br{(E-g)V_Q}\ge\dim V_Q-(2l-\ep l)=2l(k+1)-(2l-\ep l)=2kl+\ep l$. Множество $Q$ содержит $(k+1)l$ весов с учётом кратностей, следовательно,
$\bn{(E-A)Q}\le(k+1)l$, и
\eqn{\label{dimm}
\dim\br{(E-g)V_Q}-\bn{(E-A)Q}\ge2kl+\ep l-(k+1)l=(k-1+\ep)l.}

Число $k-1+\ep$ может быть отрицательным только при условии $k=\ep=0$ или, что равносильно, $A\la=\la$. Таким образом, если $A\la\ne\la$, то автоматически $\la\ne0$, и
\eqn{\label{rkQ}
\dim\br{(E-g)V_Q}\ge\bn{(E-A)Q}.}
Если же $A\la=\la$, то в $Q$ входит только вес $\la$ (с кратностью $l$), множество $(E-A)Q$ тривиально, $\bn{(E-A)Q}=0$, и неравенство~\eqref{rkQ} вновь выполняется.

Просуммировав неравенства~\eqref{rkQ} по всем независимым циклам перестановки множества изотипных компонент, получим неравенство $\rk(E-g)\ge\bn{(E-A)P}$, поскольку пространство $V$ есть прямая сумма подпространств $V_Q$, а множество $P$ есть объединение (с учётом кратностей) своих подмножеств $Q$. Следовательно, все неравенства в~\eqref{rkPg}, а также все неравенства~\eqref{rkQ} обращаются в равенства, в частности, $\bn{(E-A)P}=r+2$. Линейная оболочка $2$\д устойчивого множества $(E-A)P$ имеет размерность $r$, поэтому, согласно утверждениям~\ref{irred} и~\ref{indep}, указанное множество неразложимо, а любые его ненулевые векторы в количестве не более $r$ линейно независимы.

Как уже сказано, для каждого независимого цикла $\dim\br{(E-g)V_Q}=\bn{(E-A)Q}$. Если для какого-либо цикла $\la\ne0$, то, в силу~\eqref{dimm}, $k+\ep\le1$. Последнее неравенство верно и при $\la=0$, так как в этом случае $A\la=\la$, $k=\ep=0$. Для всякого $\la\in P^A$ изотипная компонента $V_{\la}$ переходит в себя под действием $g$ и сама образует независимый цикл. При этом множество $(E-A)Q$ тривиально, откуда $\dim\br{(E-g)V_Q}=\bn{(E-A)Q}=0$, то есть $g$ тождественно действует на $V_Q=V_{\la}$.

Теперь допустим, что включение~\eqref{Ppm} не выполняется. Тогда оператор $g$ осуществляет на множестве изотипных компонент нетождественную подстановку. Последняя имеет независимый цикл, для которого $k>0$, а поскольку $k+\ep\le1$, имеем $k=1$, $\ep=0$. Это означает, что в данный цикл входят компоненты $V_{\la^+}$ и $V_{\la^-}$, а веса $\la^{\pm}$ удовлетворяют~\eqref{Alapm}. Два ненулевых вектора $(E-A)\la^{\pm}=\la^{\pm}-\la^{\mp}$ из множества $(E-A)P$ противоположны и потому линейно зависимы. Следовательно, $r=1$, $\bn{(E-A)P}=r+2=3$, оператор $A$ является отражением, а множество $P\sm P^A$ состоит ровно из трёх весов с учётом кратностей. Веса $\la^+$ и $\la^-$ принадлежат
$P\sm P^A$, значит, $P\sm P^A=\{\la^+,\la^-,\la\}$, причём веса $\la^+$, $\la^-$ и $\la$ входят в $P$ с кратностью $1$. Все изотипные компоненты,
соответствующие весам множества $P^A$, под действием $g$ остаются на месте; $V_{\la^+}$ и $V_{\la^-}$\т меняются местами. Оставшаяся компонента $V_{\la}$ переходит в себя под действием $g$, отсюда $A\la=-\la$, и $A$\т отражение относительно гиперплоскости $\ha{\la}^{\perp}$.
\end{proof}

\begin{imp}\label{pair} Любой оператор $g\in\Om$ переводит в себя все изотипные компоненты $V_{\la}$, кроме, быть может, двух.
\end{imp}

\begin{imp}\label{Pm2} Если $\Ad(\Om)\ni-E$, то $\hn{P}=m+2$.
\end{imp}

Итак, если включение~\eqref{Ppm} не выполняется, то данному элементу $g\in\Om$ соответствует тройка весов $\{\la^+,\la^-,\la\}$, каждый из которых входит в $P$ с кратностью $1$. При этом $(E+A)(\la^+-\la^-)=(\la^+-\la^-)+(\la^--\la^+)=0$, значит,
\eqn{\label{lindep}
0\ne\la^+-\la^-\in\Ker(E+A)=\ha{\la}.}
Таким образом, между весами $\la^+$, $\la^-$ и $\la$ существует линейное соотношение, у которого все три коэффициента ненулевые. В частности, любое подпространство в $\ggt$, содержащее два вектора из тройки $\{\la^+,\la^-,\la\}$, содержит и всю её. Отметим, что
\eqn{\label{dimha}
\dim\ha{\la^+,\la^-,\la}=2,}
то есть линейная оболочка этой тройки двумерна: согласно~\eqref{Alapm}, веса $\la^{\pm}$ не лежат на прямой $\Ker(E+A)=\ha{\la}$.

\begin{theorem}\label{PA} Все неразложимые компоненты множества $P$, кроме, быть может, одной, целиком содержатся в $P^A$.
\end{theorem}

\begin{proof} Если $A=E$, то доказывать нечего.

Предположим, что $A\ne E$, а также выполнено~\eqref{Ppm}. Тогда множества $P^A$ и $P\sm P^A$ содержатся в подпространствах соответственно $\Ker(E-A)$ и $\Ker(E+A)$, имеющих тривиальное пересечение. Значит, каждое из этих множеств является объединением некоторых неразложимых компонент множества $P$. По теореме~\ref{nocent} множество $(E-A)P$ является неразложимым. Оно совпадает (с точностью до нулей) с множеством $2\br{P\sm P^A}$, поэтому множество $P\sm P^A$ также неразложимо и может включать в себя не более одной неразложимой компоненты $P$. Остальные же неразложимые компоненты полностью содержатся в $P^A$. Отметим, что если $P$ неразложимо, то (с точностью до нулевых весов) одно из множеств $P^A$ и $P\sm P^A$ пусто, а другое совпадает с $P$, а поскольку $\ha{P}=\ggt$ и $A\ne E$, имеем $A=-E$.

Если же включение~\eqref{Ppm} не выполняется, то элементу $g\in\Om$ соответствует тройка весов $\{\la^+,\la^-,\la\}\subs P$. Эти три веса лежат в одной и той же неразложимой компоненте множества $P$, так как между ними существует линейное соотношение, у которого все три коэффициента ненулевые. Все остальные неразложимые компоненты целиком содержатся в $P^A$.
\end{proof}

Вернёмся к множеству $L$ подпространств пространства $V$, определённому в п.~\ref{introd}.

\begin{lemma}\label{GVLa} Всякий оператор $g\in\Om$ действует тождественно на всех подпространствах множества $L$, кроме, быть может, одного. В частности, $gV_{\La}=V_{\La}$ для всех неразложимых компонент $\La\subs P$.
\end{lemma}

\begin{proof} При $\Ad(g)\ne E$ утверждение следует из теорем~\ref{PA} и~\ref{nocent}: все неразложимые компоненты $P$, кроме одной, лежат в $P^A$, а $g$ тождественно действует на всех $V_{\la}$, $\la\in P^A$, в том числе и на $V_0$.

Если $\Ad(g)=E$, то $\rk(E-g)=\om(g)\le2$. Оператор $g$ оставляет на месте все подпространства $V_{\la}$, действуя на каждом из них, кроме $V_0$, комплексно-линейно. Далее, имеем $V=\opluss{\la}V_{\la}$,
\equ{\label{sigdim}
\rk(E-g)=\sums{\la}\dim\br{(E-g)V_{\la}}=\sums{\la\ne0}2\dim_{\Cbb}\br{(E-g)V_{\la}}+\dim\br{(E-g)V_0}.}
Поэтому среди чисел $\dim_{\Cbb}\br{(E-g)V_{\la}}$, $\la\ne0$, и $\dim\br{(E-g)V_0}$ не может быть двух ненулевых. Следовательно, существует не более одной изотипной компоненты $V_{\la}$, на которой $g$ не действует тождественно.
\end{proof}

Теперь докажем теорему~\ref{submain}.

Предположим, что $V/G$\т гладкое многообразие. Прямое произведение подгрупп $G[W]$, $W\in L$, содержит подгруппу $G^0$, а также множество $\Om$: по лемме~\ref{GVLa} каждый элемент из $\Om$ принадлежит некоторому прямому сомножителю. Значит, согласно следствию~\ref{fst1}, указанное прямое произведение совпадает с $G$. В частности,
\eqn{\label{dirper}
G=G[W]\times G\hs{W^{\perp}}}
для любого $W\in L$.

Зафиксируем подпространство $W\in L$. В силу утверждения~\ref{MvW}, найдётся вектор $v\in W^{\perp}$, такой что $M_v\subs W$. Тогда $N_v/G_v$\т гладкое многообразие (теорема~\ref{slice}). В то же время $\ggt v\subs W^{\perp}$, $N_v\sups W$, $N_v=W\oplus\hr{N_v\cap W^{\perp}}$, причём на втором прямом слагаемом группа $G_v$ действует тождественно. Отсюда при $d:=\dim\hr{N_v\cap W^{\perp}}$ фактор $(W\oplus\R^d)/G_v$ является гладким многообразием. Кроме того, стабилизатор $G_v$ содержит подгруппу $G[W]$ и, согласно~\eqref{dirper}, совпадает с прямым произведением подгрупп $G[W]$ и $G_v\cap G\hs{W^{\perp}}$, причём последняя тождественно действует на $W$. Следовательно, $(W\oplus\R^d)/G[W]$\т гладкое многообразие.

Тем самым теорема~\ref{submain} доказана.

\begin{theorem}\label{abel} Если $\dim G=1$ и $\Ad(G)=\{E\}$, то $V/G$ не есть гладкое многообразие.
\end{theorem}

\begin{proof} Множество $P$ имеет одномерную линейную оболочку, поэтому оно неразложимо и нетривиально. Значит, $\hn{P}\ge\dim\ha{P}+2=3$, а множество $L$ может включать в себя только два подпространства: $V_0$ и $W:=V_0^{\perp}$. В силу теоремы~\ref{submain}, достаточно доказать, что ни для какого $d$ фактор $(W\oplus\R^d)/G[W]$ не является гладким многообразием. Это, в свою очередь, вытекает из теоремы~\ref{ab}: группа $G[W]$, содержащая $G^0$, одномерна, действию $G[W]\cln W$ соответствует не менее трёх ненулевых весов и ни одного нулевого, а $\Ad\br{G[W]}=\Ad(G)=\{E\}$.
\end{proof}

\begin{imp}\label{abeld} Если $\dim G=1$ и $\Ad(G)=\{E\}$, то ни при каком $d$ фактор $(V\oplus\R^d)/G$ не является гладким многообразием.
\end{imp}

\begin{proof} Представление~\eqref{act} удовлетворяет условиям теоремы~\ref{abel}, и ему соответствует $2$\д устойчивое множество весов.
\end{proof}

\section{Доказательство основной теоремы}\label{promain}

Приступим к доказательству теоремы~\ref{main}. Разобьём его на две части. В первой части мы докажем, что при $\dim G>1$ условия~\ref{m2stab}~---~\ref{adj} \textit{необходимы} для того, чтобы фактор $(V\oplus\R^d)/G$ был гладким многообразием при некотором $d$ (необходимость условия~\ref{finst} вытекает из следствия~\ref{refstab}). Во второй части будет доказано, что при $\dim G>0$ выполнения условий~\ref{m2stab}~---~\ref{finst} \textit{достаточно} для того, чтобы фактор $V/G$ был многообразием.

\subsection{Необходимость}\label{emerg}

Допустим, что множество $P$ неразложимо, $2$\д устойчиво и не содержит нулевых весов, а $m:=\dim G>1$.

Рассмотрим некоторый элемент $g\in\Om$ и положим, как и раньше, $A:=\Ad(g)$, $P^A:=P\cap\Ker(E-A)$, $r:=\rk(E-A)$. Как мы знаем, множество $(E-A)P$ неразложимо и $2$\д устойчиво, а любые его ненулевые векторы в количестве не более $r$ линейно независимы.

Если оператор $A$ не удовлетворяет~\eqref{Ppm}, то он является отражением. В этом случае $P\sm P^A=\{\la^+,\la^-,\la\}$, где $A\la=-\la$, а векторы $\la^+$ и $\la^-$ удовлетворяют соотношению~\eqref{Alapm}. Каждый из векторов $\la^+,\la^-,\la$ не пропорционален никакому другому вектору из $P$ и, в частности, входит в $P$ с кратностью $1$. Если же имеет место~\eqref{Ppm}, то $A=\pm E$, так как $P$ неразложимо. Итак, множество $\Ad(\Om)$ может
содержать только операторы $(\pm E)$ и отражения относительно гиперплоскостей в $\ggt$.

\begin{lemma}\label{tri2} Если $\hn{P}>m+2$, то всякий нетождественный оператор из $\Ad(\Om)$ является отражением относительно гиперплоскости, причём различным отражениям соответствуют тройки $\{\la^+,\la^-,\la\}$, не пересекающиеся даже с точностью до знака.
\end{lemma}

\begin{proof} Первое утверждение вытекает из следствия~\ref{Pm2}. Докажем второе утверждение.

Предположим, что двум различным отражениям $A_i=\Ad(g_i)$ ($i=1,2$, $g_i\in\Om$) соответствуют тройки $\{\la^+_i,\la^-_i,\la_i\}$, которые пересекаются с точностью до знака. Положим $P^{A_i}:=P\cap\Ker(E-A_i)$. Тогда имеем $P\sm P^{A_i}=\{\la^+_i,\la^-_i,\la_i\}$. Кроме того,
$\ha{\la_1}\ne\ha{\la_2}$, так как отражения $A_i$ относительно гиперплоскостей $\ha{\la_i}^{\perp}$ различны.

Любое подпространство в $\ggt$, содержащее два веса из тройки $\{\la^+_i,\la^-_i,\la_i\}$, содержит и всю тройку. Всякий вектор множества $P$ либо лежит в тройке $\{\la^+_1,\la^-_1,\la_1\}$, либо ортогонален вектору $\la_1$. Поэтому тройка $\{\la^+_2,\la^-_2,\la_2\}$ принадлежит объединению подпространств $\ha{\la^+_1,\la^-_1,\la_1}$ и $\ha{\la_1}^{\perp}$. По крайней мере одно из этих подпространств содержит не менее двух векторов из
тройки $\{\la^+_2,\la^-_2,\la_2\}$, а значит, и всю эту тройку. Но в тройке $\{\la^+_1,\la^-_1,\la_1\}$ ни один вектор не ортогонален $\la_1$, а тройка $\{\la^+_2,\la^-_2,\la_2\}$ с ней пересекается с точностью до знака и потому не лежит в $\ha{\la_1}^{\perp}$. Следовательно, \eqn{\label{subs21}
\{\la^+_2,\la^-_2,\la_2\}\subs\ha{\la^+_1,\la^-_1,\la_1}.}

Согласно~\eqref{dimha} и~\eqref{subs21}, тройки $\{\la^+_i,\la^-_i,\la_i\}$ имеют одну и ту же двумерную линейную оболочку $\al=\ha{\la_1,\la_2}$. В свою очередь, множество $P^{A_1}\cap P^{A_2}$, получаемое из $P$ удалением этих двух троек, лежит в $\ha{\la_1,\la_2}^{\perp}=\al^{\perp}$. В силу своей неразложимости множество $P$ совпадает с объединением троек $\{\la^+_i,\la^-_i,\la_i\}$, $i=1,2$, $\ggt=\ha{P}=\al$, $m=\dim\al=2$, $\hn{P}>m+2=4$, $\bn{P^{A_i}}=\hn{P}-3>1$. Итак, множество $P^{A_1}$ содержит более одного вектора и лежит на прямой $\ha{\la_1}^{\perp}$. Однако каждый вектор из $P$ принадлежит некоторой из троек $\{\la^+_i,\la^-_i,\la_i\}$ и не может быть пропорциональным никакому другому вектору множества $P$. Получили противоречие.
\end{proof}

\begin{lemma}\label{tri} Допустим, что всякий нетождественный оператор из $\Ad(\Om)$ является отражением относительно гиперплоскости, причём различным отражениям соответствуют тройки $\{\la^+,\la^-,\la\}$, не пересекающиеся даже с точностью до знака. Тогда ни для какого $d$ фактор $(V\oplus\R^d)/G$
не является гладким многообразием.
\end{lemma}

\begin{proof} Предположим, что $(V\oplus\R^d)/G$\т гладкое многообразие ($d\ge0$).

Пусть $\{\la_1\sco\la_p\}\subs P$\т подмножество всех весов $\la\in P$, таких что оператор отражения относительно $\ha{\la}^{\perp}$ принадлежит множеству $\Ad(\Om)$. Для каждого $i=1\sco p$ обозначим через $g_i$ элемент из $\Om$, для которого $\Ad(g_i)$\т отражение относительно гиперплоскости $\ha{\la_i}^{\perp}$ (такой существует), а через $\{\la^+_i,\la^-_i,\la_i\}$\т соответствующую ему тройку весов. Эти $p$ троек попарно не пересекаются (даже с точностью до знака).

Каждый вектор $\la_i$ ортогонален всем векторам из $P$, кроме $\la^+_i$, $\la^-_i$ и $\la_i$, в частности, $\la_i\perp\la_j$ при $i\ne j$. Векторы $\la_1\sco\la_p$ попарно ортогональны и потому линейно независимы, откуда $p\le m$.

Любой нетождественный оператор из $\Ad(\Om)$ совпадает с некоторым $\Ad(g_i)$, поэтому вся группа $\Ad(G)$ порождается отражениями $\Ad(g_i)$ (следствие~\ref{fst2}). Эти отражения попарно коммутируют, так как $\la_i$ попарно ортогональны. Следовательно, $\Ad(G)$ состоит в точности из всех операторов вида
\eqn{\label{eps1}
\br{\Ad(g_1)}^{\ep_1}\sd\br{\Ad(g_p)}^{\ep_p},}
где каждое из $\ep_i$ равно $0$ или $1$. Если для $g\in G$ представить $\Ad(g)$ в виде~\eqref{eps1}, то
\eqn{\label{eps2}
g=g_1^{\ep_1}\sd g_p^{\ep_p}\cdot g_0,}
где $\Ad(g_0)=E$. Элемент $g_i$ меняет местами $V_{\la^+_i}$ и $V_{\la^-_i}$ и оставляет на месте все остальные изотипные компоненты. Что же касается $g_0$, он вообще все изотипные компоненты переводит в себя. Значит, элемент $g$ оставляет на месте все изотипные компоненты, кроме
$V_{\la^+_1}\sco V_{\la^+_p},V_{\la^-_1}\sco V_{\la^-_p}$, а компоненты $V_{\la^+_i}$ и $V_{\la^-_i}$ переводит в себя (друг в друга), если число $\ep_i$ равно $0$ ($1$).

Допустим, что $p=m$. Тогда $\{\la_1\sco\la_m\}$\т ортогональный базис в $\ggt$, а каждый вектор $\la\in P$ ортогонален всем векторам этого базиса, кроме, быть может, одного $\la_i$ (если $\la\in\{\la^+_i,\la^-_i,\la_i\}$), отсюда $\la\in\ha{\la_i}$ для некоторого $i$. Таким образом, неразложимое множество $P$ лежит в объединении попарно ортогональных одномерных подпространств $\ha{\la_i}$, что невозможно, поскольку $\dim\ha{P}=m>1$.

В случае же $p<m$ векторы $\la^+_1\sco\la^+_p$ не порождают $\ggt$. Но все векторы из $P$, кроме $\la^-_1\sco\la^-_p$, линейно порождают $\ggt$, так как, в силу~\eqref{lindep}, вектор $\la^-_i$ линейно выражается через векторы $\la^+_i$ и $\la_i$, не принадлежащие набору $\{\la^-_1\sco \la^-_p\}$. Поэтому найдётся подмножество $Q\subs P$, содержащее все $\la^+_i$, не содержащее ни одного $\la^-_i$ и линейно порождающее гиперплоскость в $\ggt$. Рассмотрим вектор $v\in V$, который проецируется нетривиально на изотипную компоненту $V_{\la}$, если и только если $\la\in Q$. Группа $G_v$ одномерна, поскольку $\Lie G_v$\т это пересечение ядер весов множества $Q$, линейно порождающего гиперплоскость в $\ggt$. Ни один элемент $g\in G_v$ не может менять местами изотипные компоненты $V_{\la^+_i}$ и $V_{\la^-_i}$, так как вектор $v=gv$ проецируется нетривиально ровно на одну из них. Следовательно, при представлении $\Ad(g)$ в виде~\eqref{eps1} имеем $\ep_i=0$ для всех $i$, то есть $\Ad(g)=E$. Мы видим, что $\dim G_v=1$, $\Ad(G_v)=\{E\}$, отсюда представление $G_v\cln N_v$ точно, и ему (как мы уже знаем) соответствует $2$\д устойчивое множество весов. Применяя
следствие~\ref{abeld}, получаем, что фактор $(N_v\oplus\R^d)/G_v$ не является гладким многообразием. Согласно следствию~\ref{slicd}, не является им и $(V\oplus\R^d)/G$.
\end{proof}

\begin{theorem} Если $(V\oplus\R^d)/G$\т гладкое многообразие для некоторого $d$, то действие $G\cln V$ удовлетворяет условиям~\ref{m2stab}~---~\ref{adj} теоремы~\ref{main}.
\end{theorem}

\begin{proof} Множество $P$, будучи $2$\д устойчивым, содержит не менее $m+2$ весов. В силу лемм~\ref{tri2} и~\ref{tri}, оно содержит ровно $m+2$ веса, то есть условие~\ref{m2stab} выполнено. Значит, любые веса множества $P$ в количестве не более $m$ линейно независимы. В частности, кратность каждого веса равна $1$, так как $2\le m$. Таким образом, все изотипные компоненты суть одномерные комплексные пространства $V_{\la}$, $\la\in P$,
переставляемые группой $G$.

Любой элемент $g\in G$, для которого $\Ad(g)=\pm E$, переводит в себя все подпространства $V_{\la}$. Следовательно, при $\Ad(G)\ni-E$ выполняется условие~\ref{adj}, а при $\Ad(G)=\{\pm E\}$\т условия~\ref{oplu} и~\ref{adj}. Если же $\Ad(G)=\{E\}$, то $\Ad(\Om)=\{E\}$, и, согласно лемме~\ref{tri}, $(V\oplus\R^d)/G$ не есть гладкое многообразие.

Теперь предположим, что группа $\Ad(G)$ не содержится в $\{\pm E\}$. Тогда множество $\Ad(\Om)$ также не содержится в группе $\{\pm E\}$
(следствие~\ref{fst2}), и в нём есть некоторое отражение. Три вектора $\la^+$, $\la^-$, $\la$ множества~$P$, соответствующие этому отражению, линейно зависимы, откуда $m=2$, $\hn{P}=m+2=4$. Помимо этих трёх векторов, в $P$ имеется ровно один вектор $\la'$, ортогональный $\la$. Две различные прямые $\ha{\la^+}$ и $\ha{\la^-}$ на двумерной плоскости $\ggt$ симметричны относительно каждой из прямых $\ha{\la}$ и $\ha{\la'}$. Значит, среди этих четырёх весов ортогональны только $\la$ с $\la'$ и, возможно, $\la^+$ с $\la^-$. Поэтому всякое отражение в множестве $\Ad(\Om)$ переводит в себя пару прямых $\hc{\ha{\la},\ha{\la'}}$. Вообще, любой оператор из $\Ad(\Om)$ переводит в себя указанную пару прямых, так как если он не является отражением, то равен $(\pm E)$. Согласно следствию~\ref{fst2}, вся группа $\Ad(G)$ переводит в себя эту пару ортогональных прямых. Отсюда вытекает, что, во-первых, $G$ переводит в себя подпространства $V_{\la}\oplus V_{\la'}$ и $V_{\la^+}\oplus V_{\la^-}$, а во-вторых, что любой оператор в группе $\Ad(G)$ является либо поворотом на угол, кратный $\frac{\pi}{2}$, либо осевой симметрией.

Как мы видим, условие~\ref{oplu} выполнено; для проверки условия~\ref{adj} нужно доказать, что группа $\Ad(G)$ содержит оператор $(-E)$. Подгруппа поворотов в $\Ad(G)$ может включать в себя только повороты на углы, кратные $\frac{\pi}{2}$, и если она не содержит $(-E)$, то тривиальна. В таком случае $\bm{\Ad(G)}=2$, и в множестве $\Ad(\Om)$ имеется только один нетождественный оператор\т отражение. В силу леммы~\ref{tri}, фактор $(V\oplus\R^d)/G$ не является гладким многообразием, что противоречит условию.
\end{proof}

Тем самым доказана \textit{необходимость} в теореме~\ref{main}.

\subsection{Достаточность}\label{enough}

Предположим, что для действия $G\cln V$ выполнены условия~\ref{m2stab}~---~\ref{adj} теоремы~\ref{main}.

В силу условия~\ref{m2stab}, размерность $m$ группы $G$ положительна. Далее, в множестве $P$ нет нулевых весов, а все ненулевые в точности соответствуют (с учётом кратностей) подпространствам $W_i$. Последние, в свою очередь, имеют структуру одномерных комплексных пространств, на которых $G^0$ действует умножением на комплексные числа с модулем $1$. Комплексная структура (то есть умножение на $i$) переносится на пространство $V$ с прямых слагаемых $W_i$. Множество $P$ является $2$\д устойчивым и удовлетворяет равенству $\hn{P}=\dim\ha{P}+2$, поэтому оно неразложимо, любые его элементы в количестве не более $m$ линейно независимы, а в количестве не менее $m$\т линейно порождают $\ggt$. Условие~\ref{oplu} означает, что группа $G$ переставляет $W_i$, причём при $m=2$ может переставлять только $W_1$ с $W_2$ и $W_3$ с $W_4$, при $m=1$\т только $W_2$ с $W_3$, а при $m>2$ переводит в себя все $W_i$.

\begin{lemma}\label{surj} Пусть $v\in V$\т некоторый вектор. Представим его в виде $v=\suml{i=1}{m+2}w_i$, где $w_i\in W_i$. Через $I$ обозначим
множество всех чисел $i=1\sco m+2$, таких что $w_i\ne0$. Тогда
\begin{nums}{-1}
\item стабилизатор $G_v$ конечен в том и только том случае, когда $|I|\ge m$;
\item если $|I|\le m$, то $\dim G_v=m-|I|$, $\ggt v=\ha{iw_i\cln i\in I}_{\R}$, а в стабилизаторе $G_v$ найдётся элемент $g$, для которого
$\Ad(g)=-E$ и $gW_i=W_i$ при всех $i$.
\end{nums}
\end{lemma}

\begin{proof} Пусть $Q\subs P$\т множество весов, соответствующих подпространствам $W_i$, $i\in I$ (с учётом их кратностей).

Стабилизатор $G_v$ конечен, если и только если $\ha{Q}=\ggt$\т что равносильно, $|Q|\ge m$. Тем самым доказано первое утверждение. Докажем теперь второе.

Если $|I|\le m$, то $|Q|\le m$, значит, веса из $Q$ линейно независимы, гомоморфизм $\ph_Q\cln G^0\to\T^{|I|}$ сюръективен, откуда немедленно вытекают равенства $\dim G_v=m-|I|$ и $\ggt v=\ha{iw_i\cln i\in I}_{\R}$. Рассмотрим произвольный элемент $g\in G$, для которого $\Ad(g)=-E$ и $gW_i=W_i$, $i=1\sco m+2$. Поскольку $\ph_Q$\т сюръективный гомоморфизм, существует элемент $g_0\in G^0$, такой что ограничение $gg_0$ на каждое $W_i$, $i\in I$, есть осевая симметрия относительно вещественной прямой $\R w_i$. Имеем $gg_0\in G_v$, $\Ad(gg_0)=-E$ и $(gg_0)W_i=W_i$, $i=1\sco m+2$.
\end{proof}

\begin{lemma}\label{AdGE} Если $GW_i=W_i$ при всех $i=1\sco m+2$, то для любого вектора $v\in V$ с конечным стабилизатором $G_v=\ha{G_v\cap\Om}$.
\end{lemma}

\begin{proof} Так как всякий элемент $g\in G$ переводит в себя все $W_i$, оператор $\Ad(g)$ переводит каждый вес из $P$ в себя или в противоположный себе. Из неразложимости множества $P$ следует, что $\Ad(g)=\pm E$. Согласно условию~\ref{adj}, $\Ad(G)=\{\pm E\}$. Если $\Ad(g)=-E$ для некоторого $g\in G$, то $\rk\br{E-\Ad(g)}=m$. Кроме того, $g$ действует на всех $W_i$ антилинейно (осевой симметрией), отсюда $\rk(E-g)=m+2$, $\om(g)=2$. Если $\Ad(G_v)=\{\pm E\}$, то группа $G_v$ порождена всеми своими элементами $g$, для которых $\Ad(g)=-E$, а значит, $\om(g)=2$, $g\in\Om$.

Предположим теперь, что $\Ad(G_v)=\{E\}$. Пусть $I$\т множество всех $i\in\{1\sco m+2\}$, таких что $v$ нетривиально проецируется на $W_i$. Поскольку $|G_v|<\bes$ и $-E\notin\Ad(G_v)$, имеем $|I|>m$ (лемма~\ref{surj}). Для всякого элемента $g\in G_v$ оператор $\Ad(g)$ тождественный, следовательно, $g$ действует на каждом $W_i$ комплексно-линейно. В случае $i\in I$ элемент $g$ действует тождественно на комплексном пространстве $W_i$, так как сохраняет ненулевую проекцию на него вектора $v$. Поэтому ранг оператора $E-g$ вдвое больше, чем количество чисел $i\in\{1\sco m+2\}$, таких что
$g$ действует нетождественно на $W_i$, которое, в свою очередь, не превосходит $m+2-|I|\le1$. Итак, $\rk(E-g)\in\{0;2\}$ и $E-\Ad(g)=0$, откуда $g\in\Om$.
\end{proof}

\begin{lemma}\label{infst} Рассмотрим некоторый вектор $v'\in V$ с бесконечным стабилизатором. Представим его в виде $v'=\suml{i=1}{m+2}w_i$, где $w_i\in W_i$. Пусть $I$\т множество всех чисел $i=1\sco m+2$, таких что $w_i\ne0$, а $V\"$\т прямая сумма всех $W_i$ по $i\notin I$. Тогда
\begin{nums}{-1}
\item пространство $N_{v'}$ есть прямая сумма $G_{v'}$\д инвариантных подпространств $V\"$ и $\ha{w_i\cln i\in I}_{\R}$, причём последнее совпадает с
$i\ggt v'$, и на нём $G_{v'}$ действует тождественно;
\item если действие $G\cln V$ удовлетворяет всем четырём условиям~\ref{m2stab}~---~\ref{finst} теоремы~\ref{main}, то аналогичным условиям
удовлетворяет и действие $G_{v'}\cln V\"$ (вообще говоря, с другим значением $m$).
\end{nums}
\end{lemma}

\begin{proof} Если $v'=0$, то доказывать нечего. В дальнейшем будем считать, что $v'\ne0$. Обозначим через $V'$ прямую сумму всех $W_i$ по $i\in I$. Тогда $|I|>0$, $V=V'\oplus V\"$.

Пользуясь леммой~\ref{surj} и условием $|G_{v'}|=\bes$, получаем, что $0<|I|<m$, $\dim G_{v'}=m-|I|$, $\ggt v'=\ha{iw_i\cln i\in I}_{\R}$,
$i\ggt v'=\ha{w_i\cln i\in I}_{\R}$, $N_{v'}=V\"\oplus i\ggt v'$.

Всякий элемент $g\in G_{v'}$ переставляет подпространства $W_i$, $i\in I$\т те, на которые $v'$ проецируется нетривиально. Если указанные подпространства переставляются нетождественно, то $|I|\ge2$ и, в силу условия~\ref{oplu}, $m=\dim G\le2$, то есть $m=\dim G\le2\le|I|<m$, что невозможно. Следовательно, любой элемент $g\in G_{v'}$ переводит в себя все подпространства $W_i$, $i\in I$, а значит, и проекции на них вектора $v'$\т а именно, все векторы $w_i$, $i\in I$. Итак, $G_{v'}$ действует тождественно на подпространстве $\ha{w_i\cln i\in I}_{\R}$; кроме того,
$G_{v'}V'=V'$ и $G_{v'}V\"=V\"$.

Пространство $V\"$ есть прямая сумма $m+2-|I|=\dim G_{v'}+2$ одномерных комплексных пространств $W_i$, $i\notin I$. Так как группа $G_{v'}^0$ действует на $V'$ тождественно, представление $G_{v'}\cln V\"$ имеет конечное ядро неэффективности, и ему, как и представлению $G_{v'}\cln V$, соответствует $2$\д устойчивое множество весов. В стабилизаторе $G_{v'}$ найдётся элемент $g$, переводящий в себя все $W_i$ и удовлетворяющий равенству $\Ad(g)=-E$ (вновь по лемме~\ref{surj}).

Для группы $G$ выполнено по крайней мере одно из двух условий $\dim G\le2$ и $\fa i\;GW_i=W_i$, значит, соответствующее условие справедливо и для любой её подгруппы, в том числе $G_{v'}$. Далее, группа $G$ может переставлять подпространства $W_i$ только внутри фиксированных непересекающихся пар, и их (пар) количество не превосходит двух; это свойство также переносится на действие $G_{v'}\cln V\"$.

Тем самым доказано, что действие $G_{v'}\cln V\"$ удовлетворяет условиям~\ref{m2stab}~---~\ref{adj}. Осталось доказать, что выполняется условие~\ref{finst}.

Пусть $v\"\in V\"$\т произвольный вектор, причём его стабилизатор $G':=G_{v'}\cap G_{v\"}$ в группе $G_{v'}$ конечен. Нужно доказать, что группа $G'$ порождена элементами, действующими на факторпространстве $V\"/(\ggt_{v'}v\")$ псевдоотражением либо тождественно. Положим $v:=v'+v\"$.

Очевидно, что $G'\subs G_v$. Кроме того, $|G_v|<\bes$, поскольку группа $G_v^0$, переводящая в себя подпространства $V'$ и $V\"$, содержится в конечной группе $G'$. Пусть $\Om_v:=G_v\cap\Om$\т множество всех элементов группы $G_v$, действующих на факторпространстве $V/(\ggt v)$ псевдоотражением либо тождественно. По условию $\ha{\Om_v}=G_v$.

Так как группа $G'$ содержится в $G_{v'}$ и $G_v$, она действует на $i\ggt v'$ тождественно и переводит в себя подпространства $V'$, $V\"$, $\ggt v'$, $\ggt v$, а также подпространство $\ggt_{v'}v\"$. Покажем, что $\ggt_{v'}v\"\subs V\"$, $\ggt v\subs\ggt v'+V\"$, причём факторпредставления $V\"/(\ggt_{v'}v\")$ и $\br{\ggt v'+V\"}/(\ggt v)$ группы $G'$ изоморфны. Действительно, имеем $\ggt V'\subs V'$, $\ggt V\"\subs V\"$, откуда \equ{\label{intsec}
V\"\sups\ggt v\cap V\"=\ggt(v'+v\")\cap V\"=\ggt_{v'}(v'+v\")=\ggt_{v'}v\",}
\equ{\label{summa}
\ggt v\subs\ggt v+V\"=\ggt(v'+v\")+V\"=\ggt v'+V\"=\ggt v'\oplus V\",}
а факторпредставления $V\"/(\ggt v\cap V\")$ и $\br{\ggt v+V\"}/(\ggt v)$ группы $G'$ изоморфны. Группа $G'$ действует тождественно на подпространстве $i\ggt v'=\br{\ggt v'\oplus V\"}^{\perp}$. Поэтому $G'\cap\Om_v$\т это в точности множество всех элементов из $G'$, действующих на
$\br{\ggt v'\oplus V\"}/\ggt v$\т что равносильно, на $V\"/(\ggt_{v'}v\")$\т псевдоотражением либо тождественно. Таким образом, требуется доказать, что множество $G'\cap\Om_v$ порождает группу $G'$.

Группа $G_{v'}$ тождественно действует на подпространстве $\ha{w_i\cln i\in I}_{\R}$ и переводит в себя все $W_i$, $i\in I$, а переставлять может только $W_i$, $i\notin I$. Значит, $G_{v'}$ может переставлять подпространства $W_i$ только внутри фиксированной пары: если бы этих пар было две, то $m+2-|I|\ge4$, $m+2>4$, $m>2$, но при $m>2$ группа $G$ вообще не переставляет $W_i$. Так как $m+2-|I|>2$, можно считать, что эта пара включает в себя подпространства $W_{m+1}$ и $W_{m+2}$, $m+1,m+2\notin I$. Пусть $\rh$\т гомоморфизм из $G$ в циклическую группу $\Z_2$, ядро которого состоит в точности из всех элементов, переводящих в себя $W_{m+1}$ и $W_{m+2}$. Тогда $G_{v'}\cap\Ker\rh$\т подгруппа в $G_{v'}$ конечного индекса, она
переводит в себя все $W_i$, $i=1\sco m+2$, и её действие на $V\"$, как и действие $G_{v'}\cln V\"$, удовлетворяет условиям~\ref{m2stab}~---~\ref{adj}. Применим к действию $G_{v'}\cap\Ker\rh$ на $V\"$ лемму~\ref{AdGE} и получим, что конечная группа $G'\cap\Ker\rh=(G_{v'}\cap\Ker\rh)\cap G_{v\"}$
порождена элементами, действующими на $V\"/(\ggt_{v'}v\")$ псевдоотражением либо тождественно. Поэтому
\eqn{\label{Kerrh}
\ha{G'\cap\Om_v}\sups(G'\cap\Ker\rh).}

В множестве $\Om_v$ ни один элемент не может переставлять подпространства $W_i$ сразу в двух парах. В самом деле, если элемент $g\in\Om_v$ переставляет их сразу в двух (непересекающихся) парах, то $m+2\ge4$, $m\ge2$, и среди весов множества $P$ любые два линейно независимы, в частности, кратность каждого из них равна $1$, а $W_i$\т изотипные компоненты представления $G^0\cln V$. Согласно следствию~\ref{pair}, оператор $g\in\Om_v\subs\Om$ переводит в себя все изотипные компоненты $W_i$, кроме, быть может, двух.

Произвольный элемент $g\in\Om_v$, не лежащий в $\Ker\rh$, принадлежит $G'$. Действительно, он переставляет подпространства $W_{m+1}$ и $W_{m+2}$ и потому переводит в себя все остальные $W_i$. В частности, $gW_i=W_i$ при всех $i\in I$, так как $m+1,m+2\notin I$. Отсюда элемент $g\in G_v$ переводит в себя $V'$ и $V\"$, а значит, векторы $v'$ и $v\"$.

Итак, $\Om_v\subs(G'\cup\Ker\rh)$, поэтому множество $\rh(G'\cap\Om_v)$ совпадает с $\rh(\Om_v)$ с точностью до тривиального элемента, а поскольку $\ha{\Om_v}=G_v$, имеем
\equ{\label{rhG'}
\rh\br{\ha{G'\cap\Om_v}}=\ba{\rh(G'\cap\Om_v)}=\ba{\rh(\Om_v)}=\rh\br{\ha{\Om_v}}=\rh(G_v)\sups\rh(G').}
Из этого, а также из~\eqref{Kerrh} следует нужное равенство $\ha{G'\cap\Om_v}=G'$.
\end{proof}

\begin{lemma}\label{pitri} Пространство $V$ разлагается в прямую сумму $G$\д инвариантных подпространств $V'$ и $V\"$, таких что фактор единичной сферы в $V\"$ по действию $G$ гомеоморфен замкнутому шару, а все гомотопические группы фактора $\br{V'\sm\{0\}}/G$ тривиальны.
\end{lemma}

\begin{proof} Пусть $V':=W_{m+1}\oplus W_{m+2}$, $V\":=W_1\sop W_m$, тогда $V=V'\oplus V\"$, $GV'=V'$, $GV\"=V\"$. Подпространству $V\"$ соответствуют $m$ линейно независимых весов из $P$, кроме того, $G$ переставляет подпространства $W_1\sco W_m$, а при $m>2$ переводит в себя каждое из них. Согласно теореме~\ref{stor}, фактор единичной сферы в $V\"$ по действию $G$ гомеоморфен замкнутому $(m-1)$\д мерному шару.

Осталось доказать, что все гомотопические группы фактора $\br{V'\sm\{0\}}/G$ тривиальны. При $m\ge2$ это следует из теоремы~\ref{stor}: два веса из $P$, соответствующие $V'$, линейно независимы, а $G$ переставляет $W_{m+1}$ и $W_{m+2}$. В случае же $m=1$ нужное утверждение вытекает из теоремы~\ref{dtor}. В самом деле, действию одномерной группы $G$ на четырёхмерном пространстве $V'$ соответствуют два ненулевых веса. При этом существует элемент $g\in G$, переводящий в себя все $W_i$, для которого $\Ad(g)=-E$. Он действует на каждом одномерном комплексном пространстве $W_i$
антилинейно, то есть осевой симметрией, и поэтому инволютивен.
\end{proof}

\begin{lemma}\label{fstt} Если действие $G\cln V$ удовлетворяет условию~\ref{finst} теоремы~\ref{main}, то $V/G$\т многообразие.
\end{lemma}

\begin{proof} Будем доказывать утверждение индукцией по (положительной) размерности $m$ группы $G$. Допустим, что $m=\dim G>0$, а для всех групп положительных размерностей, меньших $m$, утверждение уже доказано.

Достаточно доказать, что фактор единичной сферы $S\subs V$ по действию $G$ является многообразием: останется воспользоваться леммами~\ref{pitriv} и~\ref{pitri}. Для всякого вектора $v\in S$ с конечным стабилизатором имеет место равенство $G_v=\ha{G_v\cap\Om}$, и, в силу следствия~\ref{refstab}, фактор $S/G$ локально в точке $\pi(v)$ является многообразием. Покажем, что для любого вектора $v\in S$ с бесконечным стабилизатором фактор пространства $N_v(S)=N_v\cap(T_vS)=N_v\cap\ha{v}^{\perp}$ по действию $G_v$ есть многообразие. По теореме о слайсе из этого будет следовать, что фактор $S/G$ локально в каждой точке является многообразием.

Пусть $v\in S$, $|G_v|=\bes$, $I$\т множество чисел $i\in\{1\sco m+2\}$, таких что проекция $v$ на $W_i$ ненулевая, $V\"$\т прямая сумма всех $W_i$ по $i\notin I$. Тогда $0<|I|<m$. Из леммы~\ref{infst} вытекает, что $N_v=V\"\oplus i\ggt v$, $v\in i\ggt v\subs N_v^{G_v}$, а действие $G_v\cln V\"$ удовлетворяет условиям~\ref{m2stab}~---~\ref{finst}. Отметим, что группа $G_v$ имеет положительную размерность $\dim G_v=m-|I|<m$, значит, по
предположению индукции, $V\"/G_v$\т многообразие. Кроме того, $N_v(S)\sups V\"$, $N_v(S)=V\"\oplus\br{i\ggt v\cap(T_vS)}$, причём на втором прямом слагаемом $G_v$ действует тождественно. Поэтому фактор $N_v(S)$ по действию $G_v$ является многообразием, как и фактор $V\"/G_v$.
\end{proof}

Это завершает доказательство теоремы~\ref{main}. Кроме того, из лемм~\ref{AdGE} и~\ref{fstt} легко вывести утверждение теоремы~\ref{mainp}.

\section{Случай одномерной группы}\label{1dimg}

Этот пункт будет посвящён доказательству теорем~\ref{main1}~---~\ref{GiHr}. Предположим, что $m=\dim G=1$, а множество $P$ содержит по крайней мере $3$ ненулевых веса и ни одного нулевого.

Пространство $V$ имеет комплексную структуру, стабилизатор любого его ненулевого вектора конечен, а всякое $G^0$\д инвариантное вещественное подпространство в $V$ является комплексным. Имеем $\Ad(G)\subs\{\pm E\}$, причём элемент $g\in G$, для которого оператор $\Ad(g)$ равен $E$ ($-E$), действует на $V$ линейно (антилинейно) над полем $\Cbb$. Если $(V\oplus\R^d)/G$\т гладкое многообразие для некоторого $d$, то $\Ad(G)\ne\{E\}$
(следствие~\ref{abeld}), то есть $\Ad(G)=\{\pm E\}$. Согласно следствию~\ref{fst2}, $\Ad(\Om)\ni-E$, поэтому (следствие~\ref{Pm2}) $\hn{P}=3$, $\dim_{\Cbb}V=3$. Далее будем считать, что $\dim_{\Cbb}V=3$ и $\Ad(\Om)\ni-E$, в частности, $\Ad(G)=\{\pm E\}$ и $G\ne G^0$.

Пусть $G':=\Ker\Ad\subs G$\т подгруппа всех комплексно-линейных преобразований группы $G$. Если представление $G\cln V$ неприводимо, то его ограничение на подгруппу $G'$ индекса $2$ является либо неприводимым, либо прямой суммой двух неприводимых комплексных подпредставлений равной размерности. Однако второй случай невозможен, поскольку число $\dim_{\Cbb}V=3$ нечётно. Следовательно, если представление $G\cln V$ неприводимо, то и представление $G'\cln V$ неприводимо. Группа $G$ переводит в себя все изотипные компоненты представления $G^0\cln V$. В частности, если представление $G\cln V$ неприводимо, то $V$ есть изотипная компонента, и на нём $G^0$ действует скалярными операторами $\la E$, $\la\in\T$.

Элемент $g\in G'$ принадлежит $\Om$ тогда и только тогда, когда $\rk(E-g)\in\{0;2\}$\т что равносильно, $g$ является комплексным отражением либо тождественным оператором. В свою очередь, элемент $g\in G\sm G'$ лежит в $\Om$ тогда и только тогда, когда $\rk(E-g)\in\{1;3\}$. Подпространства $\Ker(E\pm g)$ пересекаются тривиально и получаются друг из друга умножением на $i$, отсюда
\equ{\label{Lra}
g\in\Om\Lra\dim\br{\Ker(E-g)}=3\Lra\Ker(E-g)\oplus\Ker(E+g)=V\Lra g^2=E.}
Значит, $\Om\sm G'$ есть множество всех антилинейных инволюций группы $G$. Оно непусто, так как $\Ad(\Om)\ni-E$, и является объединением некоторых смежных классов $G$ по $G^0$: если $\Ad(g)=-E$, $g^2=E$, а $g_0\in G^0$, то $\Ad(gg_0)=-E$, $gg_0g=gg_0g^{-1}=g_0^{-1}$, $(gg_0)^2=(gg_0g)g_0=E$. Следовательно, $\ha{\Om}\sups G^0$. Если при этом фактор $(V\oplus\R^d)/G$ является гладким многообразием для некоторого $d$, то, в силу
следствия~\ref{fst1}, $G=\ha{G^0\cup\Om}=\ha{\Om}$. Это доказывает первое утверждение теоремы~\ref{GrHi}.

Применив к группе $G'$ предложение~\ref{fin}, получим, что все её комплексные отражения порождают в ней конечную подгруппу $H$, нормальную в $G$. По теореме~\ref{ShShT} фактор $V/H$ диффеоморфен $V$, причём любой элемент из $G'$ (из $G\sm G'$) действует на этом факторе линейно (антилинейно) над полем $\Cbb$. Представление одномерной группы $G/H$ в трёхмерном комплексном пространстве $V/H$ точно, и ему не соответствует ни одного нулевого веса. Кроме того, $\Ad(G/H)=\{\pm E\}$. Далее, отображение факторизации $V\to V/H$ может быть задано формулой $v\to\br{f_1(v),f_2(v),f_3(v)}$, где
$f_i$\т однородные многочлены, свободно порождающие алгебру инвариантов $\Cbb[V]^H$. Если $\deg f_i\ne\deg f_j$ при $i\ne j$, то группа $G$ при действии на $V/H\cong V$ переводит в себя все координатные прямые, и, согласно теореме~\ref{mainp}, $V/G$\т многообразие.

Приступим к доказательству теоремы~\ref{main1}.

\begin{lemma}\label{redref} Предположим, что группа $G$ приводима над полем $\R$, порождена множеством $\Om$ и не содержит комплексных отражений. Тогда $V/G$\т многообразие.
\end{lemma}

\begin{proof} Непустое множество $\Om\sm\{E\}$ всех антилинейных инволюций группы $G$ является объединением некоторых смежных классов $G$ по $G^0$. Пусть $\Om/G^0\subs G/G^0$\т образ этого множества при факторизации $G\to G/G^0$. Тогда если образ некоторой подгруппы в $G$ при факторизации по $G^0$ порождён своим пересечением с $\Om/G^0$, то сама подгруппа либо порождена своим пересечением с $\Om$, либо нетривиальна и содержится в $G^0$.

Поскольку группа $G$ приводима, $V$ есть прямая сумма двух $G$\д инвариантных комплексных подпространств $V'$ и $V\"$ размерностей $2$ и $1$ соответственно. Если два веса из $\ggt$, соответствующие подпространству $V'$, не равны и не противоположны, то $V$ есть прямая сумма одномерных $G$\д инвариантных комплексных подпространств $V'_1$, $V'_2$ и $V\"$, где $V'_1\oplus V'_2=V'$, и утверждение следует из теоремы~\ref{mainp}. Допустим, что
подпространству $V'$ соответствует вес $\la\cln G^0\to\T$ с кратностью $2$. Если $g\in G^0$ и $g\in\Ker\la$, то оператор $g$ линеен над $\Cbb$ и действует на $V'$ тождественно, в то время как $G$ не содержит комплексных отражений. Отсюда $\Ker\la=\{E\}$, а стабилизатор любого ненулевого вектора в $V'$ тривиально пересекается с $G^0$. Множество всех антилинейных операторов из $G$, действующих инволютивно на $V'$, совпадает с $\Om\sm\{E\}$, так как на одномерном комплексном пространстве $V\"$ всякий антилинейный ортогональный оператор инволютивен.

Докажем, что стабилизатор любого ненулевого вектора $v\in V$ порождён своим пересечением с $\Om$.

Предположим, что $v\in V'$. Пусть $\pi_0\cln V'\to V'/G^0$\т отображение факторизации. Фактор $V'/G^0$ гомеоморфен $\R^3$, конечная группа $G/G^0$ действует на этом факторе линейно, причём отражениями на нём действуют элементы множества $\Om/G^0$ и только они. Факторгруппа $G/G^0$, порождённая множеством $\Om/G^0$, действует на $V'/G^0$ как конечная группа, порождённая отражениями. Значит, стабилизатор в группе $G/G^0$ любого вектора из
$V'/G^0\cong\R^3$ порождается своими элементами, действующими на $V'/G^0$ отражениями\т что то же самое, принадлежащими $\Om/G^0$. Образ группы $G_v$ при факторизации $G$ по $G^0$ совпадает со стабилизатором в группе $G/G^0$ элемента $\pi_0(v)\in V'/G^0$ и поэтому порождается своим пересечением с $\Om/G^0$. Поскольку группа $G_v$ тривиально пересекается с $G^0$, она порождена множеством $G_v\cap\Om$.

Для ненулевого вектора $v\in V\"$ имеем $Gv=G^0v$, $G=G^0G_v$, то есть образ подгруппы $G_v\subs G$ при факторизации по $G^0$ совпадает с $G/G^0$ и порождается своим пересечением с $\Om/G^0$. Следовательно, $G_v=\ha{G_v\cap\Om}$, иначе группа $G_v$ содержалась бы в $G^0$, что невозможно, так как $G^0G_v=G\ne G^0$.

Пусть теперь вектор $v\in V$ не лежит ни в одном из подпространств $V'$ и $V\"$. Тогда $v=v'+v\"$, $v'\in V'$, $v\"\in V\"$, $v',v\"\ne0$. Всякий комплексно-линейный оператор $g\in G_v$ переводит в себя векторы $v'$ и $v\"$, линейно независимые над $\Cbb$. Отсюда подпространство $\Ker(E-g)$ имеет комплексную размерность не менее $2$, и, ввиду отсутствия в $G$ отражений, $g=E$. Итак, $G_v\cap G'=\{E\}$, $|G_v|\le2$, значит, в группе $G_v$
любой нетождественный оператор антилинеен и инволютивен, то есть принадлежит $\Om\sm\{E\}$.

Согласно следствию~\ref{refstab}, фактор единичной сферы $S\subs V$ по действию $G$ является многообразием: нами установлено, что стабилизатор произвольного вектора $v\in S$ конечен и порождён своим пересечением с $\Om$. Как уже было сказано, группа $G/G^0$ действует на линейном пространстве $V'/G^0$ как конечная группа, порождённая отражениями. В таком случае фактор $V'/G\cong(V'/G^0)/(G/G^0)$ гомеоморфен замкнутому неограниченному выпуклому многограннику, причём образ нуля при факторизации является граничной точкой. В частности, все гомотопические группы фактора $\br{V'\sm\{0\}}/G$ тривиальны. Единичная сфера в пространстве $V\"$ представляет собой в точности одну орбиту группы $G$, а её фактор по действию $G$ состоит из одной точки. В силу леммы~\ref{pitriv}, $V/G$\т многообразие.
\end{proof}

\begin{imp}\label{red} Если представление $(G/H)\cln(V/H)$ приводимо, а $G=\ha{\Om}$, то $V/G$\т многообразие.
\end{imp}

\begin{proof} Поскольку $V/G\cong(V/H)/(G/H)$, достаточно доказать, что для действия $(G/H)\cln(V/H)$ выполнены условия леммы~\ref{redref}. В самом деле, это представление приводимо, а в группе $G/H$ ни один элемент не действует на $V/H$ комплексным отражением. Наконец, условие $G=\ha{\Om}$ переносится и на группу $G/H$. Действительно, при гомоморфизме факторизации $G\to G/H$ образ подмножества $\Om\subs G$ содержится в аналогичном подмножестве группы $G/H$: антилинейные инволюции переходят вновь в антилинейные инволюции, а все остальные операторы из $\Om$ принадлежат $H$ и
переходят в тождественный оператор.
\end{proof}

Таким образом, доказано второе утверждение теоремы~\ref{main1}. Будем доказывать первое.

Допустим, что группа $G^0$ действует на $V$ скалярными операторами. Тогда естественно рассмотреть в группе $G'$ подгруппу $K:=G'\cap\SL_{\Cbb}(V)$ всех её операторов с комплексным определителем $1$. Легко видеть, что $K\nl G$, $|K|<\bes$ и $G'=G^0K$. Кроме того, отсутствие комплексных отражений в $G$ равносильно тому, что в группе $K$ любой оператор над полем $\Cbb$ либо является скалярным, либо имеет простой спектр. Если $K$\т
группа чётного порядка, то в ней есть (нетождественный) оператор порядка $2$, все его собственные значения равны $(\pm1)$ и дают в произведении $1$. Такой оператор не является скалярным и не обладает простым спектром.

Следовательно, если $G$ не содержит комплексных отражений, а $G^0$ действует на $V$ скалярными операторами, то группа $K$ имеет конечный нечётный порядок, и $G'=G^0K$. Поэтому на любой системе импримитивности представления $G\cln V$ группа $G'$ осуществляет только чётные подстановки. Если при этом группа $G'$ неприводима над $\Cbb$, то неприводима и группа $K$, пересекающая все смежные классы $G'$ по $G^0$. Все конечные неприводимые подгруппы в $\SL_3(\Cbb)$ описаны в~\cite[~гл.~V]{Blich}. Так, все примитивные группы имеют чётный порядок. Значит, группа $K$ импримитивна. На своей системе импримитивности, состоящей из трёх прямых, она осуществляет всевозможные чётные подстановки. Если такая система импримитивности не единственна, то $K=G'(3,3,3)$.

\begin{lemma}\label{irredf} Если группа $G$ неприводима над полем $\R$ и не содержит комплексных отражений, то фактор $(V\oplus\R^d)/G$ не является гладким многообразием ни для какого $d$.
\end{lemma}

\begin{proof} Из неприводимости $G$ следует, что группа $G'$ также неприводима, а $G^0$ действует на $V$ скалярными операторами. Значит, $|K|<\bes$, $G'=G^0K$. Ввиду отсутствия в $G$ комплексных отражений, группа $K$ имеет нечётный порядок, неприводима и импримитивна, а любой её нескалярный оператор имеет простой спектр. Можно считать, что $K$ имеет в качестве системы импримитивности тройку координатных прямых, осуществляет на ней всевозможные чётные подстановки и содержит оператор циклической перестановки координат. В группе $K$ все диагональные операторы образуют подгруппу индекса $3$, которая вместе с любым своим оператором $\diag(\la_1,\la_2,\la_3)$ содержит $\diag(\la_2,\la_3,\la_1)$ и $\diag(\la_3,\la_1,\la_2)$. Все
операторы группы $G$ переставляют системы импримитивности её нормальной подгруппы $K$.

Предположим, что $(V\oplus\R^d)/G$\т гладкое многообразие для некоторого $d$. Докажем, что $K=G'(3,3,3)$. В самом деле, если группа $K$ не совпадает с $G'(3,3,3)$, то имеет ровно одну систему импримитивности, состоящую из трёх прямых. Тогда все операторы из $G$ переставляют координатные прямые
$\Cbb e_i$, $1\le i\le3$. Согласно следствию~\ref{fst1}, найдётся оператор $g\in\Om\sm\{E\}$, переставляющий эти прямые нетождественно. Будучи антилинейной инволюцией, он осуществляет на их множестве транспозицию и меняет местами, не умаляя общности, прямые $\Cbb e_2$ и $\Cbb e_3$. Группа $K$ вместе с любым своим оператором $g_1=\diag(\la_1,\la_2,\la_3)$ содержит и $g_2:=gg_1g^{-1}=\diag(\ol{\la_1},\ol{\la_3},\ol{\la_2})$, а также
$g_1g_2^{-1}=\diag(\la_1^2,\la_2\la_3,\la_2\la_3)$. Оператор $g_1g_2^{-1}$ не обладает простым спектром и потому скалярный: $\la_1^2=\la_2\la_3$, $\la_1^3=\la_1\la_2\la_3=\det g_1=1$. Но операторы $\diag(\la_2,\la_3,\la_1)$ и $\diag(\la_3,\la_1,\la_2)$ тоже принадлежат $K$, следовательно, $\la_i^3=1$, $1\le i\le3$, то есть $g_1^3=E$. Итак, любой диагональный оператор в $K$ даёт при возведении в куб тождественный оператор и имеет определитель $1$. Значит, $K=G'(3,3,3)$.

Группа $K=G'(3,3,3)$ имеет порядок $27$. Её коммутант совпадает с центром $\Zc(K)$, состоит в точности из трёх скалярных операторов с определителем $1$ и изоморфен $\Z_3\cong\F_3^+$\т аддитивной группе поля $\F_3$. Факторгруппа $K/\Zc(K)$ коммутативна и изоморфна $\Z_3\times\Z_3$\т аддитивной группе двумерного пространства $\F_3^2$. Далее, отображению $K\times K\to\Zc(K)$, $(a,b)\to aba^{-1}b^{-1}$ соответствует невырожденная кососимметрическая билинейная форма $\wg\cln\F_3^2\times\F_3^2\to\F_3$. Любой автоморфизм группы $K$ действует на её коммутанте $\Zc(K)\cong\F_3^+$ умножением на число $\ep\in\F_3^*$, а на факторгруппе $K/\Zc(K)$ как на пространстве $\F_3^2$\т невырожденным линейным преобразованием, умножающим форму $\wg$ на это же число $\ep\in\F_3^*$.

Системами импримитивности представления $K\cln V$ являются ровно четыре попарно не пересекающиеся тройки прямых. Эти четыре системы находятся в биекции с четырьмя прямыми пространства $\F_3^2$. Именно, в группе $K$ подгруппа всех операторов, переводящих в себя каждую прямую некоторой (фиксированной) системы импримитивности, имеет индекс $3$, содержит центр $\Zc(K)$ и при факторизации по нему даёт прямую в $\F_3^2$. Если $g\in G$ и $\Ad(g)=\ep E$, $\ep=\pm1$, то сопряжение при помощи $g$ переводит в себя группу $K$. Полученный автоморфизм группы $K$ действует на её центре, изоморфном $\F_3^+$, умножением на число $\ep\in\F_3^*$, а соответствующее ему линейное преобразование в $\F_3^2$ умножает форму $\wg$ на $\ep\in\F_3^*$ и осуществляет на четырёхэлементном множестве прямых в $\F_3^2$ подстановку знака $\ep$. Эту же подстановку знака $\ep$ осуществляет
оператор $g$ на четырёхэлементном множестве систем импримитивности представления $K\cln V$.

Зафиксируем одну из систем импримитивности и выберем в ней прямую $\Cbb v$, $v\ne0$. Стабилизатор вектора $v$ в группе $K$ состоит из трёх элементов. Отсюда конечная группа $G_v$ нетривиальна, а множество $G_v\cap\Om$, в силу следствия~\ref{refstab}, содержит нетождественный оператор $g\in\Om\sm\{E\}$. Имеем $g(\Cbb v)=\Cbb v$, поэтому данная система импримитивности сохраняется оператором $g\in G\sm G'$, подгруппами $K$ и $G^0$, а значит, и всей группой $G$. Итак, все операторы из $G$ переводят в себя каждую из четырёх систем импримитивности, то есть осуществляют на их множестве тождественную подстановку, знак которой равен $1$. Следовательно, $\Ad(G)=\{E\}$, что неверно.
\end{proof}

\begin{imp}\label{irredh} Если представление $(G/H)\cln(V/H)$ неприводимо, то фактор $(V\oplus\R^d)/G$ не является гладким многообразием ни для какого $d$.
\end{imp}

Теперь теорема~\ref{main1} доказана полностью.

\begin{lemma}\label{redred} Если либо группа $G$ приводима над полем $\R$, либо группа $H$ неприводима над полем $\Cbb$, то представление $(G/H)\cln(V/H)$ приводимо.
\end{lemma}

\begin{proof} Вначале предположим, что $G$\т приводимая группа. Тогда $V$ есть прямая сумма ненулевых $G$\д инвариантных комплексных подпространств $V'$ и $V\"$, а каждое комплексное отражение в $G$ действует на одном из них тождественно, на другом\т отражением. В таком случае $V/H\cong(V'/H)\times(V\"/H)$, $V'/H\cong V'$, $V\"/H\cong V\"$, причём группа $G/H$ действует линейно на факторах $V'/H$ и $V\"/H$. Линейное представление группы $G/H$ в пространстве $V/H$ приводимо как прямая сумма её представлений в ненулевых пространствах $V'/H$ и $V\"/H$.

Допустим теперь, что группа $H$ неприводима. В~\cite{ShT} приведена классификация всех конечных неприводимых подгрупп $H\subs\GL_3(\Cbb)$, порождённых отражениями, и для каждой из них указан набор степеней однородных образующих алгебры инвариантов $\Cbb[V]^H$. Так, если группа $H$ примитивна, то степени однородных образующих $\Cbb[V]^H$ попарно различны. Если же $H$ импримитивна, то $H=G(pq,p,3)$ для некоторых $p,q\in\N$. Алгебра $\Cbb[V]^H$ имеет однородные образующие степеней $pq,2pq,3q$ (отметим, что при $p\ne3$ числа $pq,2pq,3q$ \textit{попарно} различны). В любом случае среди степеней однородных образующих $\Cbb[V]^H$ по крайней мере две различны, поэтому некоторая координатная прямая комплексного пространства $V/H$ инвариантна относительно $G/H$.
\end{proof}

Это завершает доказательство теоремы~\ref{GrHi}: первое утверждение было доказано в начале этого пункта, а второе следует из леммы~\ref{redred} и теоремы~\ref{main1}. Кроме того, если группа $H$ неприводима и не совпадает ни с одной из групп $G(3q,3,3)$, $q\in\N$, то $V/G$\т многообразие, так как степени однородных образующих алгебры инвариантов $\Cbb[V]^H$ попарно различны. Тем самым доказана и теорема~\ref{HG3}.

Теперь докажем теорему~\ref{GiHr}.

Предположим, что группа $H$ приводима над $\Cbb$, а $G$\т неприводима над $\R$. Тогда группа $G^0$ действует на $V$ скалярными операторами.

Допустим, что $(V\oplus\R^d)/G$\т гладкое многообразие для некоторого $d\ge0$. По теореме~\ref{main1} представление $(G/H)\cln(V/H)$ приводимо, в отличие от представления $G\cln V$, в частности, группа $H$ нетривиальна.

Представление $H\cln V$ не может быть изотипичным, иначе $V$ есть прямая сумма трёх $H$\д инвариантных прямых, группа $H$ действует на $V$ скалярными операторами, но в то же время порождена отражениями и нетривиальна. Следовательно~(\cite[~\S~9]{Serr}), представление $G\cln V$ имеет систему импримитивности, состоящую из трёх $H$\д инвариантных комплексных прямых (можно считать, что координатных), и \textit{индуцировано} представлением собственной подгруппы $\wt{G}\ne G$, $H\subs\wt{G}\subs G$, на прямой $\Cbb e_1\subs V$. Поэтому группа $H$ состоит в точности из всех диагональных
операторов $h$, таких что $h^q=E$, $q>1$, то есть $H\nl G'(q,1,3)$\т подгруппа индекса $3$. Алгебра инвариантов $\Cbb[V]^H$ имеет однородные образующие $z_i^q$ одной и той же степени $q$, и при отображении факторизации $V\to V/H$, $(z_1,z_2,z_3)\to(z_1^q,z_2^q,z_3^q)$ группа $G^0$ действует на $V/H$ скалярными операторами.

Приводимое представление $(G/H)\cln(V/H)$ индуцировано одномерным представлением подгруппы $\wt{G}/H$ на прямой $\Cbb e_1\subs V/H$. Согласно \textit{критерию неприводимости Макки}~(\cite[~\S~7]{Serr}), любой оператор из $G'$, переводящий в себя все три координатные прямые пространств $V$ и $V/H$, действует на $V/H$ скалярно. Кроме того, группа $G'$ осуществляет на множестве координатных прямых в $V$ и $V/H$ только чётные подстановки, так как ни один элемент из $G/H$ не действует на $V/H$ комплексным отражением. Можно считать, что $G'$ порождена подгруппой $G^0H$ и оператором циклической перестановки координат, а прямая $\Cbb v$, где $v:=(1,1,1)\in V/H$, инвариантна относительно $G/H$. Таким образом, $G'=G^0G'(q,1,3)=G'(q,3)\subs G(q,3)$.

Остаётся доказать, что оператор $g_{(2,3)}$ принадлежит $G$. Группа $H$ переводит в себя все три координатные прямые, отсюда, в силу следствия~\ref{fst1}, некоторая антилинейная инволюция $g\in\Om$ переставляет их нетождественно, а значит, транспозицией. Можно считать, что
$g(z_1,z_2,z_3)=(\la_1\ol{z_1},\la\ol{z_3},\la\ol{z_2})$, $|\la|=|\la_1|=1$. Поскольку вектор $v=(1,1,1)\in V/H$ под действием $g$ переходит в вектор $(\la_1^q,\la^q,\la^q)\in V/H$, имеем $\la_1^q=\la^q$, и $g_{(2,3)}\in gG^0H\subs G$.

Обратно, пусть $G=G(q,3)$. Группа $H$ состоит в точности из всех диагональных операторов $h$, таких что $h^q=E$, и отображение факторизации $V\to V/H$ можно задать формулой $(z_1,z_2,z_3)\to(z_1^q,z_2^q,z_3^q)$. Легко видеть, что представление группы $G$ в пространстве $V/H$ приводимо: комплексная
прямая, натянутая на вектор $(1,1,1)\in V/H$, остаётся инвариантной. Осталось доказать, что $\ha{\Om}=G$. Антилинейные инволюции $g_{(2,3)}$ и
$g_{(1,2)}\cln(z_1,z_2,z_3)\to(\ol{z_2},\ol{z_1},\ol{z_3})$ принадлежат $G$, поэтому группа $\ha{\Om}$ осуществляет всевозможные подстановки на множестве трёх прямых $\Cbb e_i\subs V$. Подгруппа всех операторов группы $G$, переводящих в себя все координатные прямые, совпадает с $G^0H$. Имеем $\ha{\Om}\sups H$, а также $\ha{\Om}\sups G^0$, так как $\Ad(\Om)\ni-E$, значит, $\ha{\Om}=G$.

Теорема~\ref{GiHr} доказана.

\section{<<Упрощение>> представлений}\label{simpl}

\subsection{Общая конструкция}\label{gensimp}

Предположим, что имеется абстрактная группа $G$, а также подгруппа $H\subs G$, действующая на некотором (абстрактном) множестве $F$. Через $G/H$ будем вновь обозначать множество левых смежных классов $G$ по $H$. На множестве $G\times F$ рассмотрим действие группы $H$ по правилу $h(g,f):=(gh^{-1},hf)$ и действие группы $G$ левыми сдвигами первого аргумента: $g'(g,f):=(g'g,f)$. Действие группы $G$ перестановочно с действием $H$, поскольку левые сдвиги коммутируют с правыми. Поэтому группа $G$ канонически действует на факторе $(G\times F)/H$\т то есть однородном расслоении $G\us{H}*F$\т так,
что сюръективное отображение факторизации $\pi_0\cln G\times F\to G\us{H}*F$ коммутирует с действием $G$ на множествах $G\times F$ и $G\us{H}*F$. Сюръекция же $\pi'\cln G\times F\to G/H,\,(g,f)\to gH$ постоянна на слоях $\pi_0$ и перестановочна с действием $G$ на $G\times F$ и $G/H$ левыми сдвигами. Значит, существует сюръективное отображение $\pi\cln G\us{H}*F\to G/H$, перестановочное с действием $G$ и такое что $\pi'=\pi\circ\pi_0$.

Если $f,f'\in F$\т различные элементы, то $\pi_0(e,f)$ и $\pi_0(e,f')$\т различные элементы в $G\us{H}*F$, поэтому отождествим элемент
$\pi_0(e,f)\in G\us{H}*F$ ($f\in F$) с $f$. Тогда $\pi^{-1}(eH)=\pi_0(H\times F)=\pi_0\br{\{e\}\times F}=F$, а поскольку $\pi$ перестановочно с действием $G$, имеем $\pi^{-1}(gH)=g\pi^{-1}(eH)=gF$. Итак, $\pi$\т расслоение со слоем $F$. Рассмотрим множество $S_H(G,F)$ его сечений\т отображений $s\cln G/H\to G\us{H}*F$, для которых $\pi\circ s=\id_{G/H}$. На $S_H(G,F)$ канонически действует группа $G$ по правилу $g(s):=g\circ s\circ g^{-1}$, так как
\equ{\label{pigs}
\pi\hr{gsg^{-1}}=g\pi sg^{-1}=g\circ\id_{G/H}\circ g^{-1}=\id_{G/H}.}

Если $F$\т линейное пространство и $H$ действует на $F$ линейно, то слой расслоения $\pi$\т линейное пространство. Поэтому все сечения расслоения $\pi$ также образуют векторное пространство $S_H(G,F)$, причём действие $G\cln S_H(G,F)$ линейно. Если при этом $\dim F<\bes$ и $|G/H|<\bes$, то пространство $S_H(G,F)$ конечномерно.

Ниже будет описан один способ свед\'ения исходного линейного представления $G\cln V$ к представлениям в пространствах меньшей размерности. А именно, допустим, что пространство $V$ содержит прямую сумму $S$ попарно ортогональных подпространств, переставляемых группой $G$, причём последняя действует на их (конечном) множестве $M$ транзитивно. Тогда произведение всех подгрупп вида $G[W]\subs G$ ($W\in M$) прямое. Положим $T:=S^{\perp}\subs V$. Ясно, что $GS=S$, $GT=T$, а множество подпространств $M\cup\{T\}$ является системой импримитивности представления $G\cln V$.

Теперь зафиксируем некоторое подпространство $W\in M$. Пусть $G(W):=\{g\in G\cln gW=W\}\subs G$\т стабилизатор элемента $W\in M$ в группе $G$, это подгруппа Ли конечного индекса $|G/G(W)|=|M|$, откуда $G^0\subs G(W)$. Очевидно, что $G[W]\subs G(W)$, более того, $G[W]\nl G(W)$.

Рассмотрим некоторую подгруппу $H\subs G[W]$, нормальную в $G(W)$. Она тождественно действует на $W^{\perp}$ и потому вкладывается в $\Or(W)$. Для любого $g\in G$ имеем $gHg^{-1}\subs G[gW]$. При этом если $g_1W=g_2W$, $g_i\in G$, то $g_0:=g_1^{-1}g_2\in G(W)$, $g_0Hg_0^{-1}=H$, откуда $g_2Hg_2^{-1}=g_1Hg_1^{-1}$. Итак, в группе $G$ различные подгруппы вида $gHg^{-1}$ лежат в различных группах вида $G[gW]$, $g\in G$. Значит, среди подгрупп $gHg^{-1}$ лишь конечное число различных (не более $|M|$), и их произведение $H'$ прямое. При сопряжении элементами из $G$ подгруппы вида $gHg^{-1}$ переставляются, следовательно, $H'\nl G$.

\begin{prop}\label{simp} Предположим, что фактор $W/H$ диффеоморфен векторному пространству $W'$, на котором группа $G(W)$ действует линейно. В этом случае фактор $V/H'$ диффеоморфен некоторому векторному пространству $V'$, причём $G$ действует на нём линейно.
\end{prop}

\begin{proof} Подгруппа $G(W)\subs G$ конечного индекса $|M|$ действует линейно в пространстве $W'$, значит, множество $S':=S_{G(W)}(G,W')$ является
конечномерным векторным пространством, и $G$ действует на нём линейно. Согласно условию, фактор $S/H'$ диффеоморфен конечномерному векторному пространству $S'$, и каноническое действие группы $G$ на этом факторе линейно. А так как $H'$ действует на $T$ тождественно, фактор $V/H'$ диффеоморфен конечномерному пространству $V':=S'\oplus T$, причём $G$ действует на нём линейно.
\end{proof}

Пусть $G'\subs\GL(V')$\т образ группы $G$ при полученном линейном представлении $G\cln V'$ (в условиях предложения~\ref{simp}). В таком случае фактор $V/G$ является (гладким) многообразием тогда и только тогда, когда $V'/G$\т (гладкое) многообразие. Тем самым вопрос о том, является ли фактор представления $G\cln V$ (гладким) многообразием, сводится к аналогичному вопросу для точного линейного представления $G'\cln V'$. Если при этом
$\dim H>0$, то $\dim H'>0$, $\dim V'=\dim(V/H')<\dim V$, то есть исходное представление удалось свести к представлению строго меньшей размерности.

\subsection{Переход к $2$\д устойчивому множеству весов}\label{tostab}

Вернёмся к действию группы $G$ в пространстве $V$, весам $\la\cln G^0\to\T$ и $\la\cln\ggt\to\R$, соответствующим разложению на неприводимые компоненты (с учётом кратностей). Множество этих весов в пространстве $\ggt$ вновь обозначим через $P$. Как было сказано, множество $P$ не зависит от выбора разложения (с точностью до замены знака некоторых его элементов). Будем считать его $1$\д устойчивым, иначе, согласно предложению~\ref{1st}, $V/G$\т не многообразие. Предположим, что $P$ не является $2$\д устойчивым множеством. Наша цель\т свести изучение фактора исходного представления к изучению фактора другого точного представления $G'\cln V'$, для которого $\dim V'<\dim V$. Основным инструментом послужит описанный в п.~\ref{gensimp} способ перехода к другому представлению. Для доказательства того, что какое-либо подпространство $W$ и какая-либо подгруппа $H\subs G[W]$, нормальная в $G(W)$, удовлетворяют условиям предложения~\ref{simp}, будем пользоваться примерами из п.~\ref{torum}.

Любой вес $\la\in\ggt$ можно сменить на противоположный ему, поэтому можно добиться того, чтобы все коэффициенты вида $c_{\la}$, $\la\in P$, были положительны (см. п.~\ref{vect}).

Всякий оператор $\Ad(g)$, $g\in G$, с точностью до знака переводит в себя множество $P$ и переставляет его классы эквивалентности.

Поскольку множество $P$ не является $2$\д устойчивым, в нём найдётся класс эквивалентности $N$, состоящий более чем из одного элемента. Этот класс не содержит тривиальных весов и, кроме того, вместе с любым своим весом включает в себя и все его копии из $P$. Положим $W:=V_N$. Среди образов множества $N$ под действием всех операторов $\Ad(g)$, $g\in G$, любые два либо совпадают, либо не пересекаются (с точностью до знака). Поэтому среди подпространств $gW$, $g\in G$, любые два либо совпадают, либо ортогональны. Их сумма прямая, и они транзитивно переставляются под действием $G$. Пусть $H:=G[W]\cap G^0\subs G[W]$. Эта подгруппа нормальна в $G(W)$, так как $G[W]\nl G(W)$ и $G^0\nl G$.

\begin{prop}\label{simpst} При данном выборе подпространства $W$ и подгруппы $H\subs G[W]$ фактор $W/H$ диффеоморфен некоторому векторному пространству $W'$, а группа $G(W)$ действует на нём линейно.
\end{prop}

\begin{proof} Легко видеть, что подгруппа $H=G[W]\cap G^0$ совпадает с $\caps{\la\in P\sm N}\Ker\la$ ($\la\cln G^0\to\T$\т веса), а подалгебра $\hgt:=\Lie H=\ggt[W]$\т с $\caps{\la\in P\sm N}\Ker\la$ ($\la$\т веса из $\ggt$).

Группа ограничений операторов из $G^0$ на подпространство $W$ действует на каждой его неприводимой компоненте умножением на комплексное число, значит, она естественным образом вкладывается в тор $\T^{\hn{N}}$, причём образ при вложении есть не что иное как $\ph_N(G^0)$. Что касается группы $H$, то она тоже вкладывается в $\T^{\hn{N}}$, образ при вложении\т $\ph_N(H)$.

Ограничим гомоморфизм $\ph_N$ на подгруппу Ли $H\subs G^0$. Если некоторый характер
\equ{\label{char}
\chi\cln\T^{\hn{N}}\to\T,\,z\to\prods{\la\in N}z_{\la}^{d_{\la}},\;d_{\la}\in\Z}
тривиален на всей подгруппе $\ph_N(H)\subs\T^{\hn{N}}$, то линейная функция
\equ{\label{dchar}
d\chi\cln\R^{\hn{N}}\to\R,\,x\to\sums{\la\in N}d_{\la}x_{\la}}
тривиальна на алгебре $\Lie\ph_N(H)=d\ph_N(\hgt)\subs\R^{\hn{N}}$. Последнее эквивалентно тому, что функция $\sums{\la\in N}d_{\la}\la$ равна нулю на подпространстве $\hgt=\caps{\la\in P\sm N}\Ker\la$ или, что равносильно, лежит в $\ha{P\sm N}$. Поскольку $N$\т класс эквивалентности, среди таких линейных комбинаций нетривиальная ровно одна с точностью до пропорциональности, причём её коэффициенты $c_{\la}$ ($\la\in N$) положительны. Значит, $\Lie\ph_N(H)$\т гиперплоскость в $\R^{\hn{N}}$, подгруппа $\ph_N(H)\subs\T^{\hn{N}}$ имеет коразмерность $1$ и задаётся одним нетривиальным характером с показателями степеней, пропорциональными $c_{\la}$ и поэтому не равными нулю и имеющими один и тот же знак (можно считать, что натуральными).

Пространство $W$ есть прямая сумма $\hn{N}$ попарно ортогональных неприводимых компонент, каждая из них является одномерным комплексным пространством, на котором группа $G^0$ (а значит, и $H$) действует умножением на комплексные числа. Единственная с точностью до пропорциональности нетривиальная линейная комбинация весов $\la\in N$, лежащая в линейной оболочке остальных весов из $P$, имеет коэффициенты $c_{\la}>0$. Как уже было сказано, $H$\т
подгруппа $\hn{N}$\д мерного тора, задаваемая одним характером с натуральными показателями степеней $a_{\la}$, которые пропорциональны $c_{\la}$: \equ{\label{phNH}
\ph_N(H)=\hc{z\in\T^{\hn{N}}\cln\prods{\la\in N}z_{\la}^{a_{\la}}=1}\subs\T^{\hn{N}};}
в частности, $\dim H=\hn{N}-1>0$. Для такой группы $H$, в силу теоремы~\ref{tor}, фактор $W/H$ диффеоморфен векторному пространству $W'$ размерности
$\hn{N}+1$. Кроме того, ограничение на $W$ любого оператора из $G(W)$\т это ортогональный оператор в $W$, лежащий в нормализаторе группы $H$, следовательно, он линейно действует на факторе $W/H$ (вновь по теореме~\ref{tor}).
\end{proof}

Согласно предложению~\ref{simp}, мы можем свести исходное действие $G\cln V$ к другому точному представлению $G'\cln V'$, такому что $V/G$\т (гладкое) многообразие, если и только если $V'/G'$\т (гладкое) многообразие. Далее, $\dim H>0$, откуда $\dim V'<\dim V$. Описанную процедуру будем продолжать до тех пор, пока не получим $2$\д устойчивое множество весов. Поскольку размерность пространства на каждом шаге строго понижается, рано или поздно процесс остановится, а значит, множество весов станет $2$\д устойчивым.

\end{document}